\newif\ifHYPER\global\HYPERtrue
\newif\ifINTERNAL\global\INTERNALfalse
\documentclass[10pt,b5paper,fleqn]{article}
\usepackage[utf8]{inputenc}
\usepackage{amsmath,amsfonts,amssymb}
\usepackage[mathscr]{eucal}
\usepackage{color}
\usepackage{graphicx}
\usepackage{amsthm}

\ifHYPER
\definecolor{ks-green}{rgb}{0.0,0.7,0.0}
\definecolor{ks-red}{rgb}{0.7,0.0,0.0}
\definecolor{ks-blue}{rgb}{0.0,0.0,0.7}
\usepackage{hyperref}
\hypersetup{colorlinks=true,citecolor=ks-green,linkcolor=ks-red}
\fi

\numberwithin{equation}{section}
\theoremstyle{plain}
\newtheorem{theorem}[equation]{Theorem}
\newtheorem{lemma}[equation]{Lemma}

\newtheorem{proposition}[equation]{Proposition}
\newtheorem{algorithm}[equation]{Algorithm}
\theoremstyle{definition}
\newtheorem{definition}[equation]{Definition}
\newtheorem{Example}[equation]{Example}
\newtheorem{Remark}[equation]{Remark}

\newenvironment{remark}{\emph{Remark.}}{}

\newenvironment{notation}{\emph{Notation.}}{}

\newcommand{\mathbox}[2]{%
  \parbox[c]{#1}{\centering #2}}

\newcommand{\lnum}[1]{\makebox[0.5\mathindent][r]{\textnormal{\footnotesize{#1}}}}
\newenvironment{algtest}{\bgroup\bfseries\upshape
  \tabbing
    \hspace*{\mathindent}\=
    \hspace*{1.5em}\=\hspace*{1.5em}\=%
    \hspace*{1.5em}\=\hspace*{1.5em}\= \kill \>}{%
  \endtabbing\egroup}

\newcommand{\block}[1]{\underline{#1}}

\newcommand{\tabstrut}{\rule[-0.5ex]{0pt}{2.8ex}}

\newcommand{\mthstrut}{\rule[-0.5ex]{0pt}{2.2ex}}

\newcommand\trp{^{\!\top}}
\newcommand\inv{^{-1}}
\newcommand\numN{\mathbb{N}}
\newcommand\numZ{\mathbb{Z}}
\newcommand\numQ{\mathbb{Q}}
\newcommand\numR{\mathbb{R}}

\newcommand{\freeALG}[2]{#1\langle #2\rangle}
\newcommand{\freeFLD}[2]{#1(\!\langle #2\rangle\!)}
\newcommand{\perm}{\Sigma}
\newcommand{\field}[1]{\mathbb{#1}}
\newcommand{\als}[1]{\mathcal{#1}}

\newcommand{\aclo}[1]{\overline{#1}}
\newcommand{\length}[1]{\vert #1 \vert}

\DeclareMathOperator{\rank}{rank}
\DeclareMathOperator{\pivot}{\#_\text{pb}}
\DeclareMathOperator{\linsp}{span}

\begin{document}
\title{Free Fractions: An Invitation\\
  to (applied) Free Fields}
\author{Konrad Schrempf%
  \footnote{Contact: math@versibilitas.at (Konrad Schrempf),
    \url{https://orcid.org/0000-0001-8509-009X},
    Universität Wien, Fakultät für Mathematik,
    Oskar-Morgenstern-Platz~1, 1090 Wien, Austria.
    }
  \hspace{0.2em}\href{https://orcid.org/0000-0001-8509-009X}{%
  \includegraphics[height=10pt]{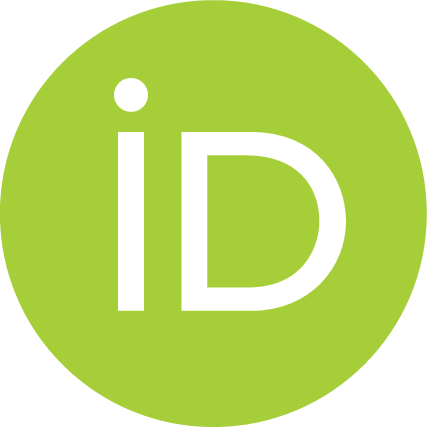}}
  }

\maketitle

\begin{abstract}
Long before we learn to construct the field of \emph{rational numbers}
(out of the ring of \emph{integers}) at university, we learn how
to calculate with \emph{fractions} at school. When it comes to numbers,
we are used to a \emph{commutative} multiplication, for example
$2\cdot 3 = 6 = 3\cdot 2$. On the other hand ---even before we can write---
we learn to talk (in a language) using \emph{words}, consisting of
purely \emph{non-commuting} letters (or symbols),
for example $xy \neq yx$ (with the concatenation as multiplication).
Now, if we combine numbers (from a field) with words (from the
\emph{free monoid} of an alphabet) we get \emph{non-commutative}
\emph{polynomials} which form a ring (with ``natural'' addition and
multiplication), namely the \emph{free associative algebra}.
Adding or multiplying polynomials is easy,
for example $(\frac{2}{3} xy + z) + \frac{1}{3} xy = xy + z$
or $2x(yx + 3z) = 2xyx + 6xz$.
Although the integers and the non-commu\-ta\-ti\-ve (nc) polynomials
look rather different, they share many properties,
for example the unique number of \emph{irreducible} factors:
$x(1-yx) = x - xyx = (1-xy)x$.
However, the construction of the \emph{universal field of fractions}
(aka ``free field'') of the \emph{free associative algebra} is
highly non-trivial (but really beautiful).
Therefore we provide techniques (building on the work of Cohn and Reutenauer)
to calculate with \emph{free fractions} (representing elements in the \emph{free field}
or ``skew field of nc rational functions'')
to be able to explore a fascinating non-commutative world.
\end{abstract}

\medskip
\emph{Keywords and 2020 Mathematics Subject Classification.}
Free associative algebra, universal field of fractions,
minimal linear representation, admissible system, rational operations,
non-commutative rational functions;
Primary 16K40, 68W30; Secondary 16S85, 16Z05

\section*{Introduction}

Since most of the literature on free fields is almost
inaccessible without a degree in mathematics and difficult
without a specialization in algebra we want to provide
an introduction with focus on the application.
One of the main hurdles is the huge number of concepts
and definitions (for precise formulations), needing a lot
of time to digest. Even if non-commutativity (as we understand
it here) is rather natural, one needs to get used to it.
Just to test ``non-commutative'' awareness:
$(x+y)^2 = \ldots$ ?

Here we restrict ourself to the simplest free fields,
coming from the embedding of the ring of non-commutative (nc)
polynomials (over a \emph{commutative} field and a \emph{finite} alphabet)
into its \emph{universal field of fractions}
\cite[Chapter~7]{Cohn2006a}% MR2246388 book
. A ``soft'' introduction is 
\cite[Section~9.3]{Cohn2003b}% MR1953987 book
. The tools (or techniques) we are going to use are mainly based
on the work of Cohn and Reutenauer \cite{Cohn1999a}% MR1723470 0218-1967
. We work \emph{directly} with (a special form of)
linear representations (aka ``free fractions'') to
\emph{add}, \emph{multiply} and \emph{invert} (non-zero)
elements in the free field. Usually one has to be 
careful and distinguish between an \emph{element} and a
\emph{representation} (of it). We know that there
are several different ``classical'' fractions for one element,
for example $\frac{-1}{2} = \frac{3}{-6} = \frac{-2}{4} = r \in \numQ$.
What we learn at school is a test to check whether two
fractions are ``equal'', that is, representing the same element.
This is the so-called \emph{word problem}.

In the general case ---for elements in the free field---,
the word problem is rather difficult because of the much
more complicated representations. Therefore it will take
some effort, to learn a number of tools from
\cite{Schrempf2017a9}% IJAC2018 0218-1967
\ (word problem, minimal inverse),
\cite{Schrempf2017b9}% JSC2018 0747-7171
\ (polynomial factorization),
\cite{Schrempf2017c9}% X171209102 arxiv
\ (general factorization theory) and
\cite{Schrempf2018a2}% X180310627 arxiv
\ (constructing minimal linear representations)
to be able to work with nc fractions.
However, these techniques enable also the implementation
in computer algebra software. For further remarks on the latter
(in German) we refer to
\cite[Section~B.5]{Schrempf2018b}% TH105 thesis
.
The perfect theoretical introduction to \emph{fractions} is
\cite{Cohn1984a}% MR758127 0024-6093
.

\medskip
Section~\ref{sec:ff.e} is meant to get acquainted with the
basic notation. 
The most important basic techniques are
presented directly in
Section~\ref{sec:ff.r} (calculating),
Section~\ref{sec:ff.f} (factorizing) and
Section~\ref{sec:ff.m} (minimizing).
In a first reading,
the (sub)sections marked with ``$\star$'' can be skipped.
Those marked with ``$\star\star$'' serve as a reference
for further reading.

\medskip
\begin{remark}
It should be noted that there are some minor differences
(in notation and definitions) between the main
publications due to the consecutive development.
The main reference (and most coherent presentation) is 
\cite{Schrempf2018b}% TH105 thesis
\ (in German), its structure of chapters and sections corresponds to
sections and subsections here.
Those who are mainly interested in polynomials
should have a look on 
\cite[Remark~1.10]{Schrempf2018a2}% X180310627 arxiv
\ before reading 
\cite{Schrempf2017b9}% JSC2018 0747-7171
. A very rich theoretical resource with focus on
free associative algebras is
\cite{Cohn1974a}% MR0379555 0021-2172
.
\end{remark}

\section{Representing Elements}\label{sec:ff.e}

First of all, we need a suitable representation of the elements
in the \emph{free field}
$\field{F} = \freeFLD{\field{K}}{X}$ of the
\emph{free associative algebra} $\freeALG{\field{K}}{X}$ over the
\emph{commutative} field $\field{K}$ (for example the rational numbers
$\numQ$ or the real numbers $\numR$) and the (finite) alphabet
$X = \{ x_1, x_2, \ldots, x_d \}$ (usually $X= \{x, y, z \}$).
Here we use a special form of a \emph{linear representation}
of Cohn and Reutenauer 
\cite{Cohn1994a}% MR1276109 0008-414X
, namely \emph{admissible linear systems}.

To illustrate such a system, we consider a \emph{linear} system
of equations $A s = v$ of dimension $n \in \numN = \{ 1, 2, 3, \ldots \}$,
that is, we have $n$ unknown components $s_1, s_2, \ldots, s_n$ in the
solution vector $s$ (and also $v$ is a column vector with $n$ rows).
If $A$ is invertible, we can write $s = A\inv v$.
Now let $n=1$ with $A = a \in \numZ \setminus \{ 0 \}$
and $v\in \numZ$ (integer entries).
Then $s = a\inv v = \frac{v}{a}$ is a representation for a rational
number $s \in \numQ$. Now, given $s_1 = \frac{v_1}{a_1}$
and $s_2 = \frac{v_2}{a_2}$, we can compute the sum $s_1 + s_2 \in \numQ$
by solving the linear system $A' s' = v'$,
\begin{displaymath}
\begin{bmatrix}
a_1 & -a_1 \\
0 & a_2
\end{bmatrix}
s'
=
\begin{bmatrix}
v_1 \\ v_2
\end{bmatrix},
\quad
s'=
\begin{bmatrix}
s_1 + s_2 \\ s_2 
\end{bmatrix}.
\end{displaymath}
(Notice the upper triangular form of the \emph{system matrix} $A'$,
the ``blocks'' in the diagonal ---here they have size $1 \times 1$---
are called \emph{pivot blocks}.)
Usually we are interested in the \emph{first} component of the
solution vector $s$. If $a_1$ and $a_2$ are invertible, then
$A$ is invertible. In that case we call $As=v$ an \emph{admissible
linear system} (ALS for short). In other words:
An ALS can represent a rational number. More general, one can
view an ALS as a ``generalized'' fraction.

Important: $A$ has to be ``invertible'' (for $n=1$ we need $A \neq 0$,
for $n >1$ we need to clarify the meaning). One can extract the
first component using the first identity (row) vector $u = e_1\trp = [1,0,\ldots,0]$,
that is, the desired element $f = s_1 = u A\inv v$.
The triple $\pi_f = (u,A,v)$ is called a \emph{linear representation}
of $f \in \field{F}$. (Recall that usually we represent a rational
number $r$ by a tuple of integers $(v,a)$,
that is, $r = \frac{v}{a}=v a\inv = a\inv v = 1 \cdot a\inv v$.
So here we could write $\pi_r = (1, a, v)$.)

\medskip
For the polynomial $f = xy+yx-yz \in \freeALG{\field{K}}{X}$
an ALS $\als{A}_f = (u,A,v)$ of dimension $n=4$ is (the
zeros are replaced by lower dots to emphasize the structure)
\begin{displaymath}
\begin{bmatrix}
1 & -x & -y & . \\
. & 1 & . & -y \\
. & . & 1 & z-x \\
. & . & . & 1
\end{bmatrix}
s =
\begin{bmatrix}
. \\ . \\ . \\ 1
\end{bmatrix},
\quad
s =
\begin{bmatrix}
xy + y(x-z) \\ y \\ x-z \\ 1
\end{bmatrix}.
\end{displaymath}
Let $A = (a_{ij})$. The solution can be easily computed
(starting from the bottom):
$s_4 = 1$ and $s_i + a_{i,i+1}s_{i+1} + \ldots + a_{i,n} s_n = 0$
for $i = 3, 2, 1$. For this special form we have invertibility
of $A$ (already over the free associative algebra $\freeALG{\field{K}}{X}$).
Is it possible to represent $f = xy+yx-yz$ by a smaller
system? And, if necessary, \emph{how} could one construct a \emph{minimal}
ALS? These are fundamental questions here, their (general) answering needs
some patience.

Later we will define the \emph{rank} of an element $f \in \field{F}$ by
the dimension of a \emph{minimal} admissible linear system (for $f$).
For a word/monomial, for example $g = xyz$, an ALS can easily be stated
(see also Proposition~\ref{pro:ff.minmon}):
\begin{displaymath}
\begin{bmatrix}
1 & -x & . & . \\
. & 1 & -y & . \\
. & . & 1 & -z \\
. & . & . & 1
\end{bmatrix}
s =
\begin{bmatrix}
. \\ . \\ . \\ 1
\end{bmatrix},
\quad
s = 
\begin{bmatrix}
xyz \\ yz \\ z \\ 1
\end{bmatrix}.
\end{displaymath}
Intuitively here it is somehow clearer (compared to the system for $f$
before) that this ALS is minimal, but we have to make that more precise.
The first goal will be to define ``simple'' \emph{rational operations}
on the level of these representations (systems), for example to
scale, to add or to multiply elements
(Proposition~\ref{pro:ff.ratop}).
That is not difficult but soon
ponderous since the systems become bigger and bigger.
And before we \emph{invert} (take the reciprocal value of) an element,
we have to ensure that this is allowed. If a system is \emph{minimal},
also that is easy. An ALS for the sum of $f_1 = 2x$ and $f_2 = 3y$ is
\begin{displaymath}
\begin{bmatrix}
1 & -x & -1 & . \\
. & 1 & . & . \\
. & . & 1 & -y \\
. & . & . & 1
\end{bmatrix}
s =
\begin{bmatrix}
. \\ 2 \\ . \\ 3
\end{bmatrix}.
\end{displaymath}
What is the solution vector $s$? Is that system minimal for $f = f_1 + f_2 = 2x +3y$?

\medskip
Let $R = \freeALG{\field{K}}{X}$. A (square) matrix
$A \in R^{n \times n}$ is called \emph{full}, if
$A = PQ$ with $P \in R^{n \times m}$ and $Q \in R^{m \times n}$
implies $m \ge n$
\cite{Cohn1999a}% MR1723470 0218-1967
. To show that the full matrices over the \emph{free associative algebra}
are those which are invertible over the free field (and vice versa)
is very difficult. For details we refer to 
\cite{Cohn2006a}% MR2246388 book
. Important for us is that we can ``address'' \emph{each} element $f$
in the free field via a \emph{linear representation}
\cite{Cohn1999a}% MR1723470 0218-1967
, that is, $\pi_f = (u,A,v)$ with (for some $n\in \numN$)
$u\trp,v \in \field{K}^{n \times 1}$, full $A \in R^{n \times n}$ with entries of
the form $\lambda_0 + \lambda_1 x_1 + \ldots + \lambda_d x_d$
with $\lambda_i \in \field{K}$ and $x_i \in X$
and $f = u A\inv v$. If $u = [1,0,\ldots,0]$ we call
$\pi_f$ an \emph{admissible linear system} and write
$\als{A}_f = \pi_f$.

\begin{remark}
The only non-invertible element in the rational numbers is zero.
In our case, the non-invertible (square) matrices are the
non-full matrices. Although the definition (of full matrices)
is simple, testing \emph{fullness} is very hard even for a
linear matrix. An example for a non-full matrix is
\begin{displaymath}
A =
\begin{bmatrix}
z & . & . \\
x & . & . \\
y & -x & 1 
\end{bmatrix}
=
\begin{bmatrix}
z & 0 \\
x & 0 \\
0 & 1
\end{bmatrix}
\begin{bmatrix}
1 & 0 & 0 \\
y & -x & 1
\end{bmatrix}.
\end{displaymath}
\end{remark}

\subsection{Free Fractions}

The main idea (of free fractions) is as simple as in the usage of ``classical''
fractions (for elements in $\numQ$):
\emph{calculating}, \emph{factorizing} and \emph{minimizing}
(or \emph{cancelling}), for example
\begin{align*}
&\frac{2}{3} \cdot \frac{3}{4}
  = \frac{6}{12} = \frac{2\cdot 3}{2\cdot 2 \cdot 3}
  = \frac{1}{2}\quad\text{or}\\
&\frac{1}{2} + \frac{3}{2}
  = \frac{4}{2}
  = \frac{2\cdot 2}{2} = 2.
\end{align*}
At some point one stops this loop and uses the fraction
(with \emph{coprime} numerator and denominator, that is,
their greatest common divisor is $1$ or $-1$).
However, the application (in our context) is not that
easy. For a concrete expression like
$f = x\inv z z\inv y z\inv = x\inv y z\inv$ one
can find a simpler (and therefore a smaller) ALS,
for example 
\begin{displaymath}
\begin{bmatrix}
x & y \\
. & z
\end{bmatrix}
s =
\begin{bmatrix}
. \\ 1
\end{bmatrix}.
\end{displaymath}
But what should one do with
$g = x - \bigl(x\inv + (y\inv - x)\inv \bigr)\inv$
from Example~\ref{ex:ff.hua}?
(Hint: $g$ is a polynomial.)

Additionally, we need \emph{minimal} admissible linear systems
for the factorization,
therefore we would run into troubles if we need the
factorization for the minimization.
The key idea to resolve this ``dependencies'' can be
guessed already in the classical setting: One can remember
the factorization of the numerator (for the product)
and the denominator (for the sum and the product).
The latter corresponds to the \emph{standard form}
(Definition~\ref{def:ff.stdals}).

\medskip
There are a lot of definitions in Section~\ref{sec:ff.rprelim}
(and even more in \cite[Section~1]{Schrempf2018a2}% X180310627 arxiv
). For an overview the mostly used will be introduced by examples. 
We take an element $f$ in the free field $\field{F}$
given by the admissible linear systems $\als{A} = (u,A,v)$
of dimension $n=4$.
(For a rational number $r \in \numQ$ we can write
$\als{A}_r = (1,a,v)$, that is, $r = 1\cdot a\inv v = \frac{v}{a}$.)
Recall that $f = u A\inv v$. If we write $s = A\inv v$,
then $f$ is the first component of the \emph{solution vector}
$s$ in the system of ``row'' equations $A s = v$.
The $n$-tuple $(s_1, s_2, \ldots, s_n)$ of entries in $s$
is called \emph{left family}. (The column solution vector
$s$ and the left family are used synonymously.)

But before we take a closer look on this system of equations,
we examine the ``column'' equations
$u = tA$, in which $f$ can be expressed as a $\field{K}$-linear combination
of the components of the row solution vector
$t = [t_1, t_2, \ldots, t_n]$.
Here, underlined entries denote \emph{static} entries,
that is, they must not be changed. If we describe (elementary)
transformations in the following, they \emph{always} refer to
the system matrix $A$.
\begin{displaymath}
\underbrace{%
\begin{bmatrix}
\underline{1} & \underline{0} & \underline{0} & \underline{0} 
\end{bmatrix}}_{\mathbox{2cm}{$u$, left hand side}}
=
\underbrace{%
\begin{bmatrix}
t_1 & t_2 & t_3 & t_4
\end{bmatrix}}_{\mathbox{2.3cm}{$t=u A\inv $, right family}}
\left.
\begin{bmatrix}
1-x & . & -x & -x \\
1 & y & 1 & -2 \\
. & 1 & . & -x \\
. & . & . & 1
\end{bmatrix}
\right\} \mathbox{2cm}{dimension $\dim\als{A} = n$}
\end{displaymath}
The equations (starting from the left) are
\begin{align*}
1 &= t_1 (1-x) + t_2, \\
0 &= t_2 y + t_3, \\
0 &= -t_1 x + t_2 \quad\text{and} \\
0 &= -t_1 x -2 t_2 -  t_3x + t_4.
\end{align*}
Instead of computing the solution $t$ immediately, we will
transform the system in such a way that this will be easier.
Now we take a look on the system $As = v$:
\begin{displaymath}
\underbrace{%
\begin{bmatrix}
1-x & . & -x & -x \\
1 & y & 1 & -2 \\
. & 1 & . & -x \\
. & . & . & 1
\end{bmatrix}}_{\mathbox{2.0cm}{$A$, system matrix}}
\underbrace{%
\begin{bmatrix}
\underline{s}_1 \\ s_2 \\ s_3 \\ s_4
\end{bmatrix}}_{\mathbox{2.0cm}{$s=A\inv v$, left family}}
=
\underbrace{%
\begin{bmatrix}
. \\ -4 \\ . \\ 2
\end{bmatrix}}_{\mathbox{1.5cm}{$v$, right hand side}}
\end{displaymath}
One equation, namely $s_4 = 2$, is especially easy to solve.
Here we have
$\kappa_1 s_1 + \kappa_2 s_2 + \kappa_3 s_3 + \kappa_4 s_4 = 1$
for $\kappa_1 = \kappa_2 = \kappa_3 = 0$ and $\kappa_4 = \frac{1}{2}$,
therefore we write $1 \in L(\als{A})$, the linear span (over $\field{K}$)
of the \emph{left family}.
(If there were not such a linear combination, we would write
$1 \not\in L(\als{A})$.)
We use an analogous notation for the linear span of the
\emph{right family} $R(\als{A})$.
Normally, we must distinguish between the element $f$ and the
representation $\als{A}$. If $\als{A}$ is \emph{minimal} (which is the
case here), we can define the \emph{rank} of $f$ as the \emph{dimension}
of $\als{A}$, $\rank f := \dim\als{A}$.
In this case we say ``$f$ is of type $(*,1)$'' or $1 \in L(f)$ if
$1 \in L(\als{A})$ respectively
``$f$ is of type $(1,*)$'' or $1 \in R(f)$ if $1 \in R(\als{A})$.

\medskip
Now we will transform this representation step by step such that
the solution of both systems of equations, that is, the computation
of $s$ and $t$, becomes easier. Those families play a crucial role
in characterizing minimality of a linear representation.
However, the goal in fact will be, that we do not have
to compute these solutions at all because, in general,
this would not help us.
Usually we write $s$ and $t$ (without its components) in
``generic'' form.
The look ``inside'' (into the representation) is only for explanation.
After the following transformation one should not forget this
``inspection'' and the computation of the ``new'' solutions $s$ and $t$
because this helps to understand the naming in \emph{left}
respectively \emph{right} family.

\begin{figure}
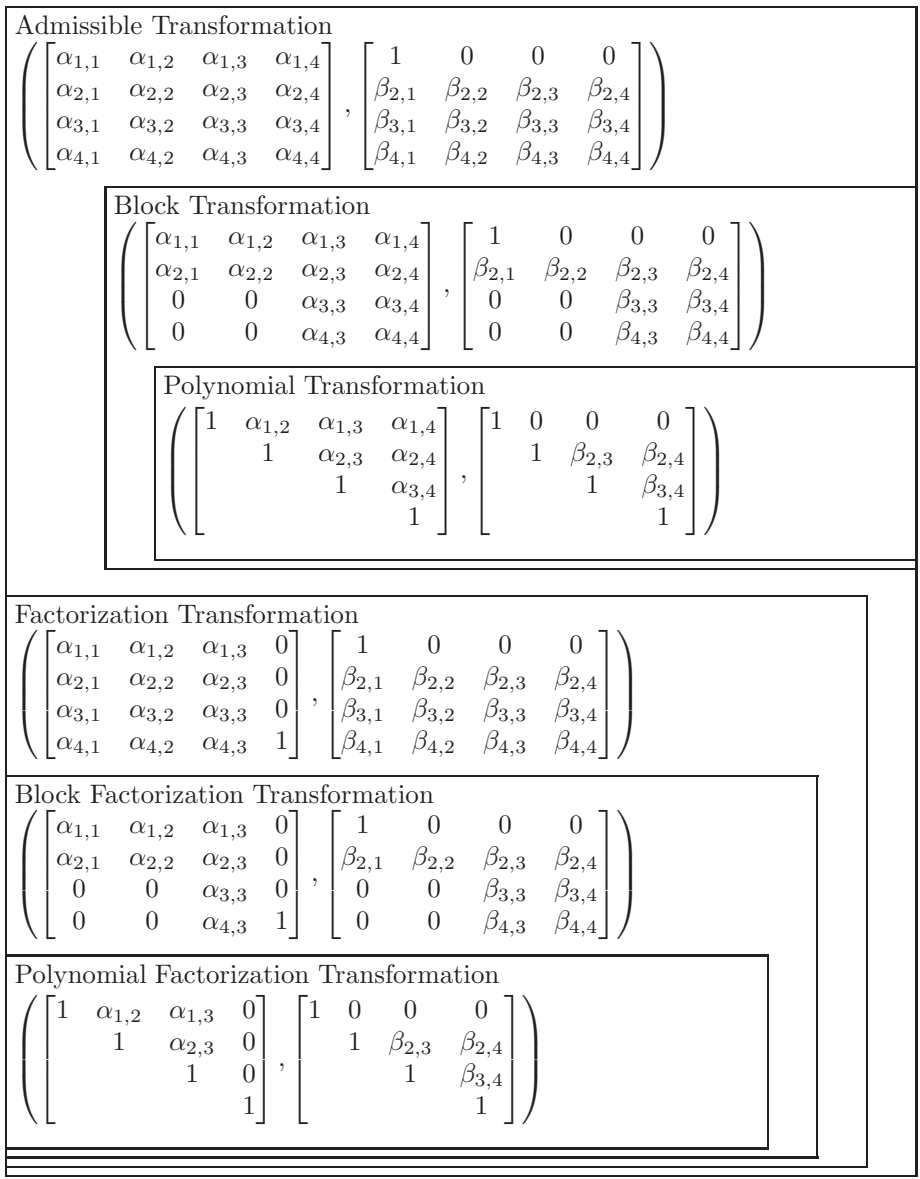

\begin{center}
\fbox{\parbox[t]{0.90\textwidth}{%
  Admissible Transformation\\[0.4ex]
  $\left(
  \begin{bmatrix}
  \alpha_{1,1} & \alpha_{1,2} & \alpha_{1,3} & \alpha_{1,4} \\
  \alpha_{2,1} & \alpha_{2,2} & \alpha_{2,3} & \alpha_{2,4} \\
  \alpha_{3,1} & \alpha_{3,2} & \alpha_{3,3} & \alpha_{3,4} \\
  \alpha_{4,1} & \alpha_{4,2} & \alpha_{4,3} & \alpha_{4,4}
  \end{bmatrix},
  \begin{bmatrix}
  1 & 0 & 0 & 0 \\
  \beta_{2,1} & \beta_{2,2} & \beta_{2,3} & \beta_{2,4} \\
  \beta_{3,1} & \beta_{3,2} & \beta_{3,3} & \beta_{3,4} \\
  \beta_{4,1} & \beta_{4,2} & \beta_{4,3} & \beta_{4,4} \\
  \end{bmatrix}
  \right)$\\[1.2ex]
\tabstrut\hspace{-\fboxsep}\hspace{-\fboxrule}\hspace{0.10\textwidth}%
\fbox{\parbox[t]{0.80\textwidth}{%
  Block Transformation\\[0.4ex]
  $\left(
  \begin{bmatrix}
  \alpha_{1,1} & \alpha_{1,2} & \alpha_{1,3} & \alpha_{1,4} \\
  \alpha_{2,1} & \alpha_{2,2} & \alpha_{2,3} & \alpha_{2,4} \\
  0 & 0 & \alpha_{3,3} & \alpha_{3,4} \\
  0 & 0 & \alpha_{4,3} & \alpha_{4,4}
  \end{bmatrix},
  \begin{bmatrix}
  1 & 0 & 0 & 0 \\
  \beta_{2,1} & \beta_{2,2} & \beta_{2,3} & \beta_{2,4} \\
  0 & 0 & \beta_{3,3} & \beta_{3,4} \\
  0 & 0 & \beta_{4,3} & \beta_{4,4} \\
  \end{bmatrix}
  \right)$\\[1.2ex]
\tabstrut\hspace{-\fboxsep}\hspace{-\fboxrule}\hspace{0.05\textwidth}%
\fbox{\parbox[t]{0.75\textwidth}{%
  Polynomial Transformation\\[0.4ex]
  $\left(
  \begin{bmatrix}
  1 & \alpha_{1,2} & \alpha_{1,3} & \alpha_{1,4} \\
    & 1 & \alpha_{2,3} & \alpha_{2,4} \\
    &   & 1 & \alpha_{3,4} \\
    &   &   & 1
  \end{bmatrix},
  \begin{bmatrix}
  1 & 0 & 0 & 0 \\
    & 1 & \beta_{2,3} & \beta_{2,4} \\
    &   & 1 & \beta_{3,4} \\
    &   &   & 1
  \end{bmatrix}
  \right)$\\[1.2ex]
}}
}}\\[2.0ex]
\tabstrut\hspace{-\fboxsep}\hspace{-\fboxrule}%
\fbox{\parbox[t]{0.85\textwidth}{%
  Factorization Transformation\\[0.4ex]
  $\left(
  \begin{bmatrix}
  \alpha_{1,1} & \alpha_{1,2} & \alpha_{1,3} & 0 \\
  \alpha_{2,1} & \alpha_{2,2} & \alpha_{2,3} & 0 \\
  \alpha_{3,1} & \alpha_{3,2} & \alpha_{3,3} & 0 \\
  \alpha_{4,1} & \alpha_{4,2} & \alpha_{4,3} & 1
  \end{bmatrix},
  \begin{bmatrix}
  1 & 0 & 0 & 0 \\
  \beta_{2,1} & \beta_{2,2} & \beta_{2,3} & \beta_{2,4} \\
  \beta_{3,1} & \beta_{3,2} & \beta_{3,3} & \beta_{3,4} \\
  \beta_{4,1} & \beta_{4,2} & \beta_{4,3} & \beta_{4,4} \\
  \end{bmatrix}
  \right)$\\[1.2ex]
\tabstrut\hspace{-\fboxsep}\hspace{-\fboxrule}%
\fbox{\parbox[t]{0.80\textwidth}{%
  Block Factorization Transformation\\[0.4ex]
  $\left(
  \begin{bmatrix}
  \alpha_{1,1} & \alpha_{1,2} & \alpha_{1,3} & 0 \\
  \alpha_{2,1} & \alpha_{2,2} & \alpha_{2,3} & 0 \\
  0 & 0 & \alpha_{3,3} & 0 \\
  0 & 0 & \alpha_{4,3} & 1
  \end{bmatrix},
  \begin{bmatrix}
  1 & 0 & 0 & 0 \\
  \beta_{2,1} & \beta_{2,2} & \beta_{2,3} & \beta_{2,4} \\
  0 & 0 & \beta_{3,3} & \beta_{3,4} \\
  0 & 0 & \beta_{4,3} & \beta_{4,4} \\
  \end{bmatrix}
  \right)$\\[1.2ex]
\tabstrut\hspace{-\fboxsep}\hspace{-\fboxrule}%
\fbox{\parbox[t]{0.75\textwidth}{%
  Polynomial Factorization Transformation\\[0.4ex]
  $\left(
  \begin{bmatrix}
  1 & \alpha_{1,2} & \alpha_{1,3} & 0 \\
    & 1 & \alpha_{2,3} & 0 \\
    &   & 1 & 0 \\
    &   &   & 1
  \end{bmatrix},
  \begin{bmatrix}
  1 & 0 & 0 & 0 \\
    & 1 & \beta_{2,3} & \beta_{2,4} \\
    &   & 1 & \beta_{3,4} \\
    &   &   & 1
  \end{bmatrix}
  \right)$\\[1.2ex]
  }}
 }}
}}
}}
\caption[Admissible Transformations (overview)]{%
The \emph{invertible} transformation matrices
$P = (\alpha_{ij}) \in \field{K}^{n \times n}$ and
$Q = (\beta_{ij}) \in \field{K}^{n \times n}$
as pair $(P,Q)$,
applied to a (not necessarily minimal)
admissible linear system $\als{A} = (u,A,v)$
of dimension $n$ for an element $f$ of the free field,
yields an \emph{equivalent} ALS $\als{A}' = P \als{A} Q = (uQ, PAQ, Pv)$
for $f$. Here $n=4$ (without loss of generality).
To ensure invertibility, we need $\det(P) \neq 0$ and
$\det(Q) \neq 0$.
}
\label{fig:ff.ueberblick}
\end{center}
\end{figure}

\medskip
Firstly we add $2$-times row~4 to row~2 (for the solution vector
$t$ this means that we subtract $2$-times $t_2$ from $t_4$).
Then we exchange columns~2 and~3 (for $s$ this means to
exchange $s_2$ and $s_3$) and subtract (the new) column~2 from
column~1.
We collect these elementary transformations in the
\emph{admissible} transformation $(P,Q)$,
that is, the first component in the solution vector $s$ does
not change, with
\begin{displaymath}
P =
\begin{bmatrix}
1 & . & . & . \\
. & 1 & . & 2 \\
. & . & 1 & . \\
. & . & . & 1
\end{bmatrix}
\quad\text{and}\quad
Q = 
\begin{bmatrix}
\underline{1} & \underline{0} & \underline{0} & \underline{0} \\
. & . & 1 & . \\
-1 & 1 & . & . \\
. & . & . & 1
\end{bmatrix}.
\end{displaymath}
(Figure~\ref{fig:ff.ueberblick} on page~\pageref{fig:ff.ueberblick}
gives an overview of different transformation matrices.)
Applying this transformation we obtain a new representation
$\als{A}' = (u',A',v')= P\als{A} Q $,
\begin{displaymath}
\als{A}' = (uQ, PAQ, Pv)
= \left(
\begin{bmatrix}
\underline{1} & \underline{0} & \underline{0} & \underline{0}
\end{bmatrix},
\begin{bmatrix}
1 & - x & . & -x \\
. & 1 & y & . \\
. & . & 1 & -x \\
. & . & . & 1
\end{bmatrix},
\begin{bmatrix}
. \\ . \\ . \\ 2
\end{bmatrix}
\right).
\end{displaymath}
The first component of the (new) solution vector $s$ is (still)
$f = 2x - 2xyx$.
Those who are not yet satisfied, can either subtract row~3 from row~1
or column~2 from column~4 and imagine our element
alternatively as $x(2-2yx)$ or $(1-xy)2x$.
This will be closer investigated in Section~\ref{sec:ff.f}.
For polynomials we always find such a form with
$n$ (scalar) ``pivot blocks'' of size $1 \times 1$.
This is not possible in general, but we will try to obtain
small pivot blocks. Either by factorization (Section~\ref{sec:ff.f})
or by ``abstract'' refinement (Section~\ref{sec:ff.m}).
But we should not worry here. The examples in the beginning
are such that we can easily minimize them by ``hand''
respectively check their minimality.

A last note concerning the system matrix $A$.
We always write it in the compact form with (at most) \emph{linear}
entries (of nc polynomials).
In fact, $A$ can also be interpreted as \emph{linear matrix pencil}
$A = (A_0, A_1, \ldots, A_d)$
with coefficient matrices $A_i \in \field{K}^{n \times n}$
for an alphabet $X = \{ x_1, \ldots, x_d \}$,
also written as $A = A_0 + A_1 x_1 + \ldots + A_d x_d$.
For an implementation one can use a list of
(square) matrices of size $n+1$.
For the example $(x-xyx)\inv$ from the beginning of Section~\ref{sec:ff.f}
with respect to the monomials $(1, x, y)$ we have
\begin{align*}
\als{A} &= (u,A,v) = \left(
\begin{bmatrix}
1 & . & .
\end{bmatrix},
\begin{bmatrix}
x & 1 & . \\
. & y & -1 \\
. & -1 & x
\end{bmatrix},
\begin{bmatrix}
. \\ . \\ 1
\end{bmatrix}
\right)\\
&\text{``$=$''}
\begin{bmatrix}
0 & u \\
v & A
\end{bmatrix}
=
\left(
\hspace{1.48em}\left|\hspace{-1.8em}%
\begin{bmatrix}
0 & 1 & 0 & 0 \\\hline\tabstrut
0 & 0 & 1 & 0 \\
0 & 0 & 0 & -1 \\
1 & 0 & -1 & 0
\end{bmatrix},
\right.
\hspace{1.48em}\left|\hspace{-1.8em}%
\begin{bmatrix}
0 & 0 & 0 & 0 \\\hline\tabstrut
0 & 1 & 0 & 0 \\
0 & 0 & 0 & 0 \\
0 & 0 & 0 & 1
\end{bmatrix},
\right.
\hspace{1.48em}\left|\hspace{-1.8em}%
\begin{bmatrix}
0 & 0 & 0 & 0 \\\hline\tabstrut
0 & 0 & 0 & 0 \\
0 & 0 & 1 & 0 \\
0 & 0 & 0 & 0
\end{bmatrix}
\right.
\right).
\end{align*}

\subsection{Left and Right Minimization Steps}

For practical computations we repeatedly have to make
admissible linear systems smaller. In concrete situations
it is possible to minimize them. Later, in Section~\ref{sec:ff.m}
we will see that there are some subtle details behind the
rather simple looking (\emph{left} and \emph{right}) ``minimization steps''.
Let us take a closer look on the example $2x + 3y$ from before:
\begin{displaymath}
\begin{bmatrix}
1 & -x & -1 & . \\
. & 1 & . & . \\
. & . & 1 & -y \\
. & . & . & 1
\end{bmatrix}
s =
\begin{bmatrix}
. \\ 2 \\ . \\ 3
\end{bmatrix},
\quad
s =
\begin{bmatrix}
2x + 3y \\ 2 \\ 3y \\ 3
\end{bmatrix}.
\end{displaymath}
First we try a ``left'' minimization step, that is,
eliminate a component of the left family. For that
we subtract $\frac{2}{3}$-times row~4 from row~2 and
add $\frac{2}{3}$-times column~2 to column~4:
\begin{displaymath}
\begin{bmatrix}
1 & -x & -1 & -\frac{2}{3} x \\
. & 1 & 0 & 0 \\
. & . & 1 & -y \\
. & . & . & 1
\end{bmatrix}
s =
\begin{bmatrix}
. \\ 0 \\ . \\ 3
\end{bmatrix},
\quad
s =
\begin{bmatrix}
2x + 3y \\ 0 \\ 3y \\ 3
\end{bmatrix}.
\end{displaymath}
The second row reads $s_2 = 0$.
That is, for the solution $s_1$ there is \emph{no contribution}
from (the new) $s_2$.
Therefore we can remove the equation $s_2=0$ and the variable
$s_2$ from our system of equations.
Hence we get the following (not yet minimal) ALS for $2x + 3y$:
\begin{displaymath}
\begin{bmatrix}
1 & -1 & -\frac{2}{3} x \\
. & 1 & -y \\
. & . & 1 
\end{bmatrix}
s =
\begin{bmatrix}
. \\ . \\ 3
\end{bmatrix},
\quad
s =
\begin{bmatrix}
2x + 3y \\ 3y \\ 3
\end{bmatrix}.
\end{displaymath}
It is obvious that now it is possible to apply a ``right''
minimization step to eliminate $t_2$ (in the right family).
In fact it is not necessary to compute the left or the right
family at all to ``minimize'' (without checking minimality).

\medskip
Minimality of a linear representation can be characterized
by $\field{K}$-linear independence of the entries of the column
solution vector $s = A\inv v$ (the \emph{left family}) and
$\field{K}$-linear independence of the entries of the row
solution vector $t = u A\inv$ (the \emph{right family})
\cite[Proposition~4.7]{Cohn1994a}% MR1276109 0008-414X
.

Since in general this is not easy to check we will
investigate conditions (on the structure of the system matrix)
in Section~\ref{sec:ff.m}
such that we can guarantee minimality if no more
(block) row and column minimization steps are possible.

\begin{Example}\label{ex:ff.blockmin}
Sometimes a minimization is only possible in ``blocks''.
Now we consider the ALS $\als{A} = (u,A,v)$%
\footnote{This ALS can be constructed in the following way:
One starts with a \emph{minimal} ALS of dimension~$3$ for the
monomial $xy$ (Proposition~\ref{pro:ff.minmon}).
Since $xy$ is of type~$(1,1)$, one can immediately ``add''
$z$ in the upper right entry of the system matrix to get
a \emph{minimal} ALS for $xy-z$.
For the inverse we use the \emph{minimal inverse}
(Theorem~\ref{thr:ff.mininv}).
And finally, using the multiplication 
(Proposition~\ref{pro:ff.ratop})
we obtain an ALS of dimension~$5$ for $f f\inv$.}
for $f f\inv = 1$
with $f = xy - z$,
\begin{displaymath}
\als{A} = \left(
\begin{bmatrix}
1 & . & . & . & . 
\end{bmatrix},
\begin{bmatrix}
1 & -x & z & . & . \\
. & 1 & -y & . & . \\
. & . & 1 & -1 & . \\
. & . & . & y & -1 \\
. & . & . & -z & x
\end{bmatrix},
\begin{bmatrix}
. \\ . \\ . \\ . \\ 1
\end{bmatrix}
\right).
\end{displaymath}
Here we can create an upper right block of zeros of size $3 \times 2$ 
in $A$ by (as a first step) adding column~3 to column~4 and
row~4 to row~2:
\begin{displaymath}
\als{A}' = \left(
\begin{bmatrix}
1 & . & . & . & . 
\end{bmatrix},
\begin{bmatrix}
1 & -x & z & z & . \\
. & 1 & -y & . & -1 \\
. & . & 1 & 0 & 0 \\
. & . & . & y & -1 \\
. & . & . & -z & x
\end{bmatrix},
\begin{bmatrix}
. \\ . \\ . \\ . \\ 1
\end{bmatrix}
\right).
\end{displaymath}
And (as a second step) adding column~2 to column~5 and
row~5 to row~1:
\begin{displaymath}
\als{A}'' = \left(
\begin{bmatrix}
1 & . & . & . & . 
\end{bmatrix},
\begin{bmatrix}
1 & -x & z & 0 & 0 \\
. & 1 & -y & 0 & 0 \\
. & . & 1 & 0 & 0 \\
. & . & . & y & -1 \\
. & . & . & -z & x
\end{bmatrix},
\begin{bmatrix}
1 \\ . \\ . \\ . \\ 1
\end{bmatrix}
\right).
\end{displaymath}
Now we can invert the lower $2 \times 2$ diagonal block
(over the free field $\field{F}$) and obtain
$t_4'' = t_5'' = 0$ (due to the zeros in the corresponding
entries in $u$).
Hence we get the (non-minimal) ALS of dimension~3,
\begin{displaymath}
\als{A}''' = \left(
\begin{bmatrix}
1 & . & . 
\end{bmatrix},
\begin{bmatrix}
1 & -x & z \\
. & 1 & -y \\
. & . & 1
\end{bmatrix},
\begin{bmatrix}
1 \\ 0 \\ 0 
\end{bmatrix}
\right).
\end{displaymath}
Notice, that the lower entries in the right hand side are zero.
Therefore a left block minimization step yields immediately
the minimal system $\als{A}''' = (1,[1],1)$ for
$1 \in \field{F}$.
\end{Example}

\medskip
\begin{remark}
The other case $f\inv f = 1$ is somewhat more difficult
because we must not change the first component in the left family.
The trick here is, to work with an ``extended'' ALS for
$1 \cdot f\inv f$, using Proposition~\ref{pro:ff.ratop}
to multiply ``$1$ from the left''.
For details and illustration see
\cite[Remark~4.3 respectively Example~4.5]{Schrempf2018a2}% X180310627 arxiv
.
\end{remark}

\begin{remark}
In some cases it is possible to do a left and a right
minimization step \emph{simultaneously}.
This is used in 
\cite[Example~5.4]{Schrempf2018a2}% X180310627 arxiv
\ to compute the left greatest common divisor of two polynomials $p$ and $q$
by minimizing an ALS for $p\inv q$.
\end{remark}

\section{Calculating}\label{sec:ff.r}

One of the main parts of this section is the
construction of a \emph{minimal} admissible linear system
for the inverse (of an element in the free field)
in Section~\ref{sec:ff.mininv}.
The following (simple) construction (of an ALS for the
inverse) is from Proposition~\ref{pro:ff.ratop}.
We assume that we have given the inverse of a monomial
$f = xyz$ by the ALS
$\als{A}' = (u',A',v')$,
\begin{displaymath}
\begin{bmatrix}
z & -1 & . \\
. & y & -1 \\
. & . & x
\end{bmatrix}
s =
\begin{bmatrix}
. \\ . \\ 1
\end{bmatrix},
\quad
s =
\begin{bmatrix}
z\inv y\inv x\inv \\
y\inv x\inv \\
x\inv
\end{bmatrix}.
\end{displaymath}
Checking also the $\field{K}$-linear independence
of the right family, the minimality is clear immediately.
A (minimal) ALS for $f$ is given by
\begin{displaymath}
\begin{bmatrix}
. & z & -1 & . \\
. & . & y & -1 \\
-1 & . & . & x \\
. & 1 & . & . 
\end{bmatrix}
s =
\begin{bmatrix}
. \\ . \\ . \\ 1
\end{bmatrix},
\quad
s =
\begin{bmatrix}
xyz \\
1 \\
z \\
yz
\end{bmatrix},
\end{displaymath}
with $-v$ in the upper \emph{left} and $u$ in the \emph{lower} right part
of the (new) system matrix. To get the form from
Proposition~\ref{pro:ff.minmon}
we are already used to, we have to reverse the rows $1,2,3$ and
columns $2,3,4$ and multiply the rows $1,2,3$ by $-1$.
As a new system $\als{A} = (u,A,v)$ for $f$ we obtain
\begin{displaymath}
\als{A} = \left(
\begin{bmatrix}
1 & . & . & . 
\end{bmatrix},
\begin{bmatrix}
1 & -x & . & . \\
. & 1 & -y & . \\
. & . & 1 & -z \\
. & . & . & 1
\end{bmatrix},
\begin{bmatrix}
. \\ . \\ . \\ 1
\end{bmatrix}
\right).
\end{displaymath}
Here it is immediate, that $1 \in L(f)$ and $1 \in R(f)$,
that is, $f$ is \emph{of type}~$(1,1)$.
The application of the inverse from
Proposition~\ref{pro:ff.ratop}
again yields an ALS for $f\inv$, however with dimension~$5$
already. Therefore an important (technical) task will be to
detect ``special'' forms (of the system matrices).

To be able to minimize, we would like to have a
very ``simple'' structure, that is, the
(diagonal) pivot blocks should be as small as possible.
For the example here, an ALS for a monomial,
this is respected by the \emph{minimal inverse}
(Theorem~\ref{thr:ff.mininv}).

\medskip
\begin{notation}
The set of the natural numbers is denoted by $\numN = \{ 1,2,\ldots \}$,
that including zero by $\numN_0$.
Zero entries in matrices are usually replaced by (lower) dots
to emphasize the structure of the non-zero entries
unless they result from transformations where there
were possibly non-zero entries before.
We denote by $I_n$ the identity matrix
and $\perm_n$ the permutation matrix that reverses the order
of rows/columns (of size $n$)
respectively $I$ and $\perm$ if the size is clear from the context.
By $v\trp$ we denote the transpose of a vector $v$.
Given an ALS $\als{A} = (u,A,v)$ with $v = [0,\ldots,0,\lambda]\trp$
we also write $\als{A} = (1,A,\lambda)$.
\end{notation}

\subsection{Preliminaries}\label{sec:ff.rprelim}

Let $\field{K}$ be a \emph{commutative} field,
$\aclo{\field{K}}$ its algebraic closure and
$X = \{ x_1, x_2, \ldots, x_d\}$ be a \emph{finite} (non-empty) alphabet.
$\freeALG{\field{K}}{X}$ denotes the \emph{free associative
algebra} (or \emph{free $\field{K}$-algebra})
and $\field{F} = \freeFLD{\field{K}}{X}$ its \emph{universal field of
fractions} (or ``free field'') 
\cite{Cohn1995a}% MR1349108 book
,
\cite{Cohn1999a}% MR1723470 0218-1967
. An element in $\freeALG{\field{K}}{X}$ is called (non-commutative or nc)
\emph{polynomial}.
In our examples the alphabet is usually $X=\{x,y,z\}$.
Including the algebra of \emph{nc rational series}
we have the following chain of inclusions:
\begin{displaymath}
\field{K}\subsetneq \freeALG{\field{K}}{X}
  \subsetneq \field{K}^{\text{rat}}\langle\!\langle X\rangle\!\rangle
  \subsetneq \freeFLD{\field{K}}{X} =: \field{F}.
\end{displaymath}
The \emph{free monoid} $X^*$ generated by $X$
is the set of all \emph{finite words}
$x_{i_1} x_{i_2} \cdots x_{i_n}$ with $i_k \in \{ 1,2,\ldots, d \}$.
An element of the alphabet is called \emph{letter},
one of the free monoid \emph{word}.
The multiplication on $X^*$ is the \emph{concatenation}
of words, that is,
$(x_{i_1} \cdots x_{i_m})\cdot (x_{j_1} \cdots x_{j_n})
= x_{i_1} \cdots x_{i_m} x_{j_1} \cdots x_{j_n}$,
with neutral element $1$, the \emph{empty word}.
The \emph{length} of a word $w=x_{i_1} x_{i_2} \cdots x_{i_m}$ is $m$,
denoted by $\length{w} = m$.
For a detailed introduction see
\cite[Chapter~1]{Berstel2011a}% MR2760561 book
.

\begin{definition}[Inner Rank, Full Matrix
\cite{Cohn2006a,Cohn1999a}% MR1723470 0218-1967
]\label{def:ff.full}
Let $R = \freeALG{\field{K}}{X}$.
Given a matrix $A \in R^{n \times n}$, the \emph{inner rank}
of $A$ is the smallest number $m\in \numN$
such that there exists a factorization
$A = T U$ with $T \in R^{n \times m}$ and
$U \in R^{m \times n}$.
The matrix $A$ is called \emph{full} if $m = n$,
\emph{non-full} otherwise.
\end{definition}

\begin{definition}[Linear Representations, Dimension, Rank
\cite{Cohn1994a,Cohn1999a}% MR1276109 0008-414X
]\label{def:ff.rep}
Let $f \in \field{F}$.
A \emph{linear representation} of $f$ is a triple $\pi_f = (u,A,v)$ with
$u\trp,v \in \field{K}^{n \times 1}$, full
$A = A_0 \otimes 1 + A_1 \otimes x_1 + \ldots
+ A_d \otimes x_d$, that is, $A$ is invertible over $\field{F}$,
$A_\ell \in \field{K}^{n\times n}$ for $\ell \in \{ 0,1,\ldots, d \}$
and $f = u A\inv v$.
The \emph{dimension} of $\pi_f$ is $\dim \, (u,A,v) = n$.
It is called \emph{minimal} if $A$ has the smallest possible dimension
among all linear representations of $f$.
The ``empty'' representation $\pi = (,,)$ is
the minimal one of $0 \in \field{F}$ with $\dim \pi = 0$.
Let $f \in \field{F}$ and $\pi$ be a \emph{minimal}
linear representation of $f$.
Then the \emph{rank} of $f$ is
defined as $\rank f = \dim \pi$.
\end{definition}

\begin{definition}[Left and Right Families
\cite{Cohn1994a}% MR1276109 0008-414X
]\label{def:ff.family}
Let $\pi=(u,A,v)$ be a linear representation of $f \in \field{F}$
of dimension $n$.
The families $( s_1, s_2, \ldots, s_n )\subseteq \field{F}$
with $s_i = (A\inv v)_i$
and $( t_1, t_2, \ldots, t_n )\subseteq \field{F}$
with $t_j = (u A\inv)_j$
are called \emph{left family} and \emph{right family} respectively.
$L(\pi) = \linsp \{ s_1, s_2, \ldots, s_n \}$ and
$R(\pi) = \linsp \{ t_1, t_2, \ldots, t_n \}$
denote their linear spans (over $\field{K}$).
\end{definition}

\begin{proposition}[%
\protect{\cite[Proposition~4.7]{Cohn1994a}% MR1276109 0008-414X
}]\label{pro:ff.cohn94.47}
A representation $\pi=(u,A,v)$ of an element $f \in \field{F}$
is minimal if and only if both, the left family
and the right family, are $\field{K}$-linearly independent.
In this case, $L(\pi)$ and $R(\pi)$ depend only on $f$.
\end{proposition}

\begin{notation}
For any two \emph{minimal} linear representations $\pi_1$ and $\pi_2$
of some element $f \in \field{F}$ we have $1 \in L(\pi_1)$
if and only if $1 \in L(\pi_2)$ because otherwise they could not
be transformed into each other by invertible matrices over $\field{K}$.
By $1 \in L(f)$ (respectively $1 \in R(f)$)
we denote $1\in L(\pi)$ (respectively $1 \in R(\pi)$)
for any \emph{minimal} $\pi$ of $f$.
\end{notation}

\begin{definition}[Element Types \protect{%
\cite[Definition~2.10]{Schrempf2017c9}% IEJA001 1306-6048
}]\label{def:ff.typele}
An element $f \in \field{F}$ is called \emph{of type} $(1,*)$
(respectively $(0,*)$) if $1 \in R(f)$ (respectively $1 \notin R(f)$).
It is called \emph{of type} $(*,1)$ (respectively $(*,0)$)
if $1 \in L(f)$ (respectively $1 \notin L(f)$).
Both subtypes can be combined.
\end{definition}

\begin{remark}
The following definition is a special case
of the more general \emph{admissible systems}
\cite[Section~7]{Cohn2006a}% MR2246388 book
\ and the slightly more general \emph{linear representations}
\cite{Cohn1994a}% MR1276109 0008-414X
.
\end{remark}

\begin{definition}[Admissible Linear Systems, Admissible Transformations
\cite{Schrempf2017a9}% IJAC2018 0218-1967
]\label{def:ff.als}
A linear representation $\als{A} = (u,A,v)$ of $f \in \field{F}$
is called \emph{admissible linear system} (ALS) for $f$,
written also as $A s = v$,
if $u=e_1=[1,0,\ldots,0]$.
The element $f$ is then the first component
of the (unique) solution vector $s$.
Given a linear representation $\als{A} = (u,A,v)$
of dimension $n$ of $f \in \field{F}$
and invertible matrices $P,Q \in \field{K}^{n\times n}$,
the transformed $P\als{A}Q = (uQ, PAQ, Pv)$ is
again a linear representation (of $f$).
If $\als{A}$ is an ALS,
the transformation $(P,Q)$ is called
\emph{admissible} if the first row of $Q$ is $e_1 = [1,0,\ldots,0]$.
\end{definition}

\begin{definition}\label{def:ff.reg}
Let $M = M_1 \otimes x_1 + \ldots + M_d \otimes x_d$
with $M_i \in \field{K}^{n \times n}$ for some $n\in \numN$.
An element in $\field{F}$ is called \emph{regular}
if it has a linear representation $(u,A,v)$ with $A = I - M$,
that is, $A_0 = I$ in Definition~\ref{def:ff.rep},
or equivalently, if $A_0$ is regular (invertible).
\end{definition}

\begin{definition}[Polynomial ALS and Transformation
\protect{\cite[Definition~24]{Schrempf2017b9}% JSC2018 0747-7171
}]\label{def:ff.psals}
An ALS $\als{A} = (u,A,v)$ of dimension $n$
with system matrix $A = (a_{ij})$
for a non-zero polynomial $0 \neq p \in \freeALG{\field{K}}{X}$ 
is called
\emph{polynomial}, if
\begin{itemize}
\item[(1)] $v = [0,\ldots,0,\lambda]\trp$ for some $\lambda \in\field{K}$ and
\item[(2)] $a_{ii}=1$ for $i=1,2,\ldots, n$ and $a_{ij}=0$ for $i>j$,
  that is, $A$ is upper triangular.
\end{itemize}
An admissible transformation $(P,Q)$ for an ALS $\als{A}$
is called \emph{polynomial} if it has the form
\begin{displaymath}
(P,Q) = \left(
\begin{bmatrix}
1 & \alpha_{1,2} & \ldots & \alpha_{1,n-1} & \alpha_{1,n} \\
  & \ddots & \ddots & \vdots & \vdots \\
  &   & 1 & \alpha_{n-2,n-1} & \alpha_{n-2,n} \\
  &   &   & 1 & \alpha_{n-1,n} \\
  &   &   &   & 1
\end{bmatrix},
\begin{bmatrix}
1 & 0 & 0 & \ldots & 0 \\
  & 1 & \beta_{2,3} & \ldots & \beta_{2,n} \\
  &   & 1 & \ddots & \vdots  \\
  &   &   & \ddots & \beta_{n-1,n} \\
  &   &   &   & 1 \\
\end{bmatrix}
\right).
\end{displaymath}
If additionally $\alpha_{1,n} = \alpha_{2,n} = \ldots = \alpha_{n-1,n} = 0$
then $(P,Q)$ is called \emph{polynomial factorization transformation}.
See also Figure~\ref{fig:ff.ueberblick} on page~\pageref{fig:ff.ueberblick}.
\end{definition}

\subsection{Minimal Systems}\label{sec:ff.minsys}

The main idea is to start with minimal admissible linear
systems and construct minimal ones for the \emph{rational operations}
(scalar multiplication, sum, product, inverse).
We already have seen the
\emph{minimal monomial}:

\begin{proposition}[Minimal Monomial
\protect{\cite[Proposition~4.1]{Schrempf2017a9}% IJAC2018 0218-1967
}]\label{pro:ff.minmon}
Let $k \in \numN$ and $f= x_{i_1} x_{i_2} \cdots x_{i_k}$ be a monomial in
$\freeALG{\field{K}}{X} \subseteq \freeFLD{\field{K}}{X}$.
Then
\begin{displaymath}
\als{A} = \left(
\begin{bmatrix}
1 & .  & \cdots & .
\end{bmatrix},
\begin{bmatrix}
1 & -x_{i_1} \\
  & 1 & -x_{i_2} \\
  & & \ddots & \ddots \\
  & & & 1 & -x_{i_k} \\
  & & & & 1
\end{bmatrix},
\begin{bmatrix}
. \\ .  \\ \vdots \\ .  \\ 1
\end{bmatrix}
\right)
\end{displaymath}
is a \emph{minimal} (polynomial) ALS of dimension $\dim\als{A} = k+1$.
\end{proposition}

More general it is possible to state \emph{minimal}
systems for a class of polynomials by a (generalized)
``companion'' system
\cite[Section~3]{Schrempf2017b9}% JSC2018 0747-7171
.

\subsection{Rational Operations}\label{sec:ff.ratop}

``Basic'' rational operations (on the level of admissible linear systems)
are easy to formulate. For the multiplication we can provide alternative
constructions yielding minimal admissible linear systems immediately
in special cases, for example the \emph{minimal polynomial multiplication}
(Proposition~\ref{pro:ff.minmul}).

\begin{proposition}[Rational Operations
\cite{Cohn1999a}% MR1723470 0218-1967
]\label{pro:ff.ratop}
Let $0\neq f,g \in \field{F}$ be given by the
admissible linear systems $\als{A}_f = (u_f, A_f, v_f)$
and $\als{A}_g = (u_g, A_g, v_g)$ respectively
and let $0\neq \mu \in \field{K}$.
Then admissible linear systems for the rational operations
can be obtained as follows:

\smallskip\noindent
The scalar multiplication
$\mu f$ is given by
\begin{displaymath}
\mu \als{A}_f =
\bigl( u_f, A_f, \mu v_f \bigr).
\end{displaymath}
The sum $f + g$ is given by
\begin{displaymath}
\als{A}_f + \als{A}_g =
\left(
\begin{bmatrix}
u_f & . 
\end{bmatrix},
\begin{bmatrix}
A_f & -A_f u_f\trp u_g \\
. & A_g
\end{bmatrix}, 
\begin{bmatrix} v_f \\ v_g \end{bmatrix}
\right).
\end{displaymath}
The product $fg$ is given by
\begin{displaymath}
\als{A}_f \cdot \als{A}_g =
\left(
\begin{bmatrix}
u_f & . 
\end{bmatrix},
\begin{bmatrix}
A_f & -v_f u_g \\
. & A_g
\end{bmatrix},
\begin{bmatrix}
. \\ v_g
\end{bmatrix}
\right).
\end{displaymath}
And the inverse $f\inv$ is given by
\begin{displaymath}
\als{A}_f\inv =
\left(
\begin{bmatrix}
1 & . 
\end{bmatrix},
\begin{bmatrix}
-v_f & A_f \\
. & u_f
\end{bmatrix},
\begin{bmatrix}
. \\ 1
\end{bmatrix}
\right).
\end{displaymath}
\end{proposition}

\begin{lemma}[for Type~$(1,*)$
\protect{\cite[Lemma~4.12]{Schrempf2017a9}% IJAC2018 0218-1967
}]\label{lem:ff.for2}
Let $\als{A} = (u,A,v)$ be a \emph{minimal} ALS
with $\dim \als{A} = n \ge 2$ and $1\in R(\als{A})$. Then there
exists an admissible transformation $(P,Q)$ such that
the first column of $PAQ$ is $[1,0,\ldots,0]\trp$ and
$Pv = [0,\ldots,0,\lambda]\trp$ for some $\lambda \in \field{K}$.
\end{lemma}

\begin{proposition}[Multiplication Type $(*,1)$ \protect{%
\cite[Proposition~3.12]{Schrempf2017c9}% IEJA001 1306-6048
}]\label{pro:ff.mul1}
Let $f,g\in \field{F} \setminus \field{K}$ be given by the
admissible linear systems $\als{A}_f = (u_f, A_f, v_f) = (1, A_f, \lambda_f)$ 
of dimension $n_f$ and
$\als{A}_g = (u_g, A_g, v_g) = (1, A_g, \lambda_g)$
of dimension $n_g$ of the form
\begin{displaymath}
\als{A}_g = \left(
\begin{bmatrix}
1 & . & .
\end{bmatrix},
\begin{bmatrix}
1 & b'& b \\
. & B & b'' \\
. & c' & c
\end{bmatrix},
\begin{bmatrix}
. \\ . \\ \lambda_g
\end{bmatrix}
\right)
\end{displaymath}
respectively.
Then an ALS for $fg$ of dimension $n = n_f + n_g -1$ is given by
\begin{displaymath}
\als{A} = \left(
\begin{bmatrix}
u_f & . & . 
\end{bmatrix},
\begin{bmatrix}
A_f & e_{n_f} \lambda_f b' & e_{n_f} \lambda_f b \\
. & B & b'' \\
. & c' & c
\end{bmatrix},
\begin{bmatrix}
. \\ . \\ \lambda_g
\end{bmatrix}
\right).
\end{displaymath}
\end{proposition}

\subsection{Disjoint Addition$^\star$}\label{sec:ff.djadd}

For \emph{disjoint} elements $f,g \in \field{F}$
\cite{Cohn1999a}% MR1723470 0218-1967
, that is, $\rank(f+g) = \rank(f) + \rank(g)$,
the addition from
Proposition~\ref{pro:ff.ratop} is minimal.
For further details we refer to the remarks after
\cite[Definition~3.2]{Schrempf2017c9}% IEJA001 1306-6048
.
An important result of Cohn and Reutenauer
is the \emph{primary decomposition} (of elements in the
free field) 
\cite[Theorem~2.3]{Cohn1999a}% MR1723470 0218-1967
.

\subsection{Minimal Inverse}\label{sec:ff.mininv}

The derivation of the \emph{minimal inverse} in 
\cite[Section~4]{Schrempf2017a9}% IJAC2018 0218-1967
\ consists of two major steps (motivated in the beginning
of this section): keeping the form for $f = (f\inv)\inv$
and distinguishing different cases to ensure minimality.
Notice especially the remark before 
\cite[Theorem~4.13]{Schrempf2017a9}% IJAC2018 0218-1967
\ how to transfer admissible linear systems into the appropriate
form.

\begin{theorem}[Minimal Inverse
\protect{\cite[Theorem~4.13]{Schrempf2017a9}% IJAC2018 0218-1967
}]\label{thr:ff.mininv}
Let $f \in \field{F} \setminus \field{K}$ be given by
the \emph{minimal} admissible linear system $\als{A} = (u, A, v)$
of dimension $n$.
Then a \emph{minimal} ALS for $f\inv$ is given
in the following way:

\smallskip\noindent
$f$ of type $(1,1)$ yields $f\inv$ of type $(0,0)$ with $\dim(\als{A}') = n-1$:
\begin{displaymath}
\als{A}'=\left(1, 
\begin{bmatrix}
-\lambda \perm b'' & -\perm B \perm \\
-\lambda b & - b'\perm
\end{bmatrix},
1 \right) \quad\text{for}\quad
\als{A} = \left(1, 
\begin{bmatrix}
1 & b' & b \\
. & B & b'' \\
. & . & 1
\end{bmatrix},
\lambda \right).
\end{displaymath}
$f$ of type $(1,0)$ yields $f\inv$ of type $(1,0)$ with $\dim(\als{A}') = n$:
\begin{displaymath}
\als{A}' = \left(1,
\begin{bmatrix}
1 & - \frac{1}{\lambda} c & - \frac{1}{\lambda} c' \perm \\
. & -\perm b'' & -\perm B \perm \\
. & -b & -b'\perm 
\end{bmatrix},
1 \right)
\quad\text{for}\quad
\als{A} = \left(1,
\begin{bmatrix}
1 & b' & b \\
. & B & b'' \\
. & c' & c
\end{bmatrix},
\lambda \right).
\end{displaymath}
$f$ of type $(0,1)$ yields $f\inv$ of type $(0,1)$ with $\dim(\als{A}') = n$:
\begin{displaymath}
\als{A}' = \left( 1,
\begin{bmatrix}
-\lambda \perm b'' & -\perm B \perm & -\perm a' \\
-\lambda b & -b'\perm & -a \\
. & . & 1 
\end{bmatrix},
1 \right)
\quad\text{for}\quad
\als{A} = \left(1,
\begin{bmatrix}
a & b' & b \\
a' & B & b'' \\
. & . & 1
\end{bmatrix},
\lambda \right).
\end{displaymath}
$f$ of type $(0,0)$ yields $f\inv$ of type $(1,1)$ with $\dim(\als{A}') = n+1$:
\begin{displaymath}
\als{A}' = \left(1,
\begin{bmatrix}
\perm v & -\perm A \perm \\
. & u \perm
\end{bmatrix},
1 \right).
\end{displaymath}
(Recall that the permutation matrix $\perm$ reverses the order of rows/columns.)
\end{theorem}

\subsection{Rational Identities}\label{sec:ff.ratid}

Using the minimal inverse (Theorem~\ref{thr:ff.mininv})
and the rational operations (Proposition~\ref{pro:ff.ratop})
one can already show non-trivial rational identities very systematically by ``hand''.
The following proof is from 
\cite[Section~5]{Schrempf2018a2}% X180310627 arxiv
.

\begin{Example}[Hua's Identity
\cite{Amitsur1966a}% MR0191912 0021-8693
]\label{ex:ff.hua}
We have:
\begin{displaymath}
x - \bigl(x\inv + (y\inv - x)\inv \bigr)\inv = xyx.
\end{displaymath}
\end{Example}

\begin{proof}
Minimal admissible linear systems for $y\inv$ and $x$ are
\begin{displaymath}
\begin{bmatrix}
y
\end{bmatrix}
s =
\begin{bmatrix}
1
\end{bmatrix}
\quad\text{and}\quad
\begin{bmatrix}
1 & -x \\
. & 1 
\end{bmatrix}
s =
\begin{bmatrix}
. \\ 1
\end{bmatrix}
\end{displaymath}
respectively.
The ALS for the difference $y\inv -x$,
\begin{displaymath}
\begin{bmatrix}
y & -y & . \\
. & 1 & -x \\
. & . & 1 
\end{bmatrix}
s = 
\begin{bmatrix}
1 \\ . \\ -1
\end{bmatrix}, \quad
s =
\begin{bmatrix}
y\inv - x \\
-x \\
-1
\end{bmatrix}, \quad
t = 
\begin{bmatrix}
y\inv & -1 & y\inv - x
\end{bmatrix},
\end{displaymath}
is minimal because the left family $s$ is $\field{K}$-linearly
independent and the right family $t$ is $\field{K}$-linearly
independent.
Clearly we have $1\in R(y\inv -x)$.
Thus, by Lemma~\ref{lem:ff.for2},
there exists an admissible transformation
\begin{displaymath}
(P,Q) = \left(
\begin{bmatrix}
. & 1 & . \\
1 & . & 1 \\
. & . & 1
\end{bmatrix},
\begin{bmatrix}
1 & . & . \\
1 & 1 & . \\
. & . & 1
\end{bmatrix}
\right),
\end{displaymath}
that yields the ALS
\begin{displaymath}
\begin{bmatrix}
1 & 1 & -x \\
. & -y & 1 \\
. & . & 1 
\end{bmatrix}
s =
\begin{bmatrix}
. \\ . \\ -1
\end{bmatrix}.
\end{displaymath}
Now we can apply the inverse of type~$(1,1)$:
\begin{displaymath}
\begin{bmatrix}
1 & y \\
-x & -1
\end{bmatrix}
s =
\begin{bmatrix}
. \\ 1
\end{bmatrix},
\quad
s = 
\begin{bmatrix}
(y\inv -x)\inv \\
-(1-xy)\inv
\end{bmatrix}.
\end{displaymath}
This system represents a regular element
$(y\inv - x)\inv = (1-yx)\inv y$,
and therefore can be transformed into a regular ALS
(Definition~\ref{def:ff.reg})
by scaling row~2 by $-1$.
Then we add $x\inv$ ``from the left'':
\begin{displaymath}
\begin{bmatrix}
x & -x & . \\
. & 1 & y \\
. & x & 1
\end{bmatrix}
s = 
\begin{bmatrix}
1 \\ . \\ -1
\end{bmatrix}, \quad
s =
\begin{bmatrix}
x\inv + (y\inv -x)\inv \\
(y\inv -x)\inv \\
-(1-xy)\inv
\end{bmatrix}.
\end{displaymath}
This system is minimal and ---after adding row~3 to row~1
(to eliminate the non-zero entry in the right hand side)---
we apply the (minimal) inverse of type~$(0,0)$:
\begin{displaymath}
\begin{bmatrix}
-1 & -1 & -x & . \\
. & -y & -1 & . \\
. & -1 & 0 & -x \\
. & . & . & 1
\end{bmatrix}
s =
\begin{bmatrix}
. \\ . \\ . \\ 1
\end{bmatrix}.
\end{displaymath}
Now we multiply row~1 and the columns~2 and~3 by $-1$
and exchange columns~2 and~3 to get the following system:
\begin{displaymath}
\begin{bmatrix}
1 & -x & -1 & . \\
. & 1 & y & . \\
. & . & 1 & -x \\
. & . & . & 1
\end{bmatrix}
s =
\begin{bmatrix}
. \\ . \\ . \\ 1
\end{bmatrix}, \quad
s = 
\begin{bmatrix}
x - xyx \\ -yx \\ x \\ 1 
\end{bmatrix}.
\end{displaymath}
The next step would be a scaling by $-1$ and
the addition of $x$ (by Proposition~\ref{pro:ff.ratop}).
With two minimization steps we would reach again minimality.
Alternatively we can add a \emph{linear} term to a polynomial
(in a polynomial ALS)
---depending on the entry $v_n$ in the right hand side---
directly in the upper right entry of the
system matrix:
\begin{displaymath}
\begin{bmatrix}
1 & -x & -1 & x \\
. & 1 & y & . \\
. & . & 1 & -x \\
. & . & . & 1
\end{bmatrix}
s =
\begin{bmatrix}
. \\ . \\ . \\ 1
\end{bmatrix}, \quad
s = 
\begin{bmatrix}
 - xyx \\ -yx \\ x \\ 1 
\end{bmatrix}.
\end{displaymath}
\end{proof}

\section{Factorizing$^\star$}\label{sec:ff.f}

Since the whole factorization theory originated from a ``small''
problem of the minimization of linear representations,
it should lead as a thread through this section.
Somehow this theory has become independent and is interesting
now from a purely algebraic point of view since it
enables to view the free field as a ``ring''.
Not in the trivial sense, where each field is a ring,
but using the richer ``structure'' by combining the
non-commutative factorization theory and the embedding
of non-commutative rings (to be more precise: \emph{free ideal rings},
FIRs
\cite{Cohn2006a}% MR2246388 book
) into their respective universal field of fraction.
There are a lot of open questions, for example,
is the free field a ``similarity unique factorization domain''?
Or, is the extension of the ``classical'' factorization theory
(in free associative algebras) to the free field ---assuming that
polynomial atoms (and their inverse) remain irreducible--- unique?

\medskip
To not loose the thread, we come back to a simple example:
Assume that we have given an element $f$ by the
admissible linear system $\als{A}_f$,
\begin{displaymath}
\begin{bmatrix}
x & 1 & . \\
. & y & -1 \\
. & -1 & x
\end{bmatrix}
s =
\begin{bmatrix}
. \\ . \\ 1
\end{bmatrix}.
\end{displaymath}
By Proposition~\ref{pro:ff.mul1}
we construct an ALS $\als{A}$ for $fx$, namely
\begin{displaymath}
\begin{bmatrix}
x & 1 & . & . \\
. & y & -1 & . \\
. & -1 & x & -x \\
. & . & . & 1
\end{bmatrix}
s =
\begin{bmatrix}
. \\ . \\ . \\ 1
\end{bmatrix}.
\end{displaymath}
Is $\als{A}$ minimal?
Now we repeat this step for $f$ given by
a different system $\als{A}_f'$ and construct
again a system $\als{A}'$ for $fx$, namely
\begin{displaymath}
\begin{bmatrix}
x & 1 & . & . \\
1 & y & -1 & . \\
. & . & x & -x \\
. & . & . & 1
\end{bmatrix}
s =
\begin{bmatrix}
. \\ . \\ . \\ 1
\end{bmatrix},
\end{displaymath}
in which one can read $\als{A}_f'$ directly in the upper left
$3 \times 3$ block of the system matrix.
Here it is immediate that row/column~3 can be eliminated
after adding column~3 to column~4.
Therefore $\als{A}'$ and hence $\als{A}$ cannot be minimal.
The connection to fac\-torization will become much clearer in
\cite[Example~30]{Schrempf2017b9}% JSC2018 0747-7171
, as soon as one verifies by the minimal inverse
that $f = (pq)\inv$ for $p=x$ and $q=1-yx$.

The (lower left) $2 \times 1$ block of zeros in the system matrix
of $\als{A}_f$ becomes an upper right block of zeros
in the system matrix of $\als{A}_f\inv$, the standard inverse of $\als{A}_f$,
\begin{displaymath}
\begin{bmatrix}
1 & -x & 1 & 0 \\
. & 1 & -y & 0 \\
. & . & -1 & -x \\
. & . & . & -1
\end{bmatrix}
s =
\begin{bmatrix}
0 \\ 0 \\ . \\ 1
\end{bmatrix},
\end{displaymath}
which is minimal here because $f$ is of type~$(0,0)$
and $\als{A}_f$ is minimal.
And this upper right block of zeros is that one
coming from multiplication $(1,*)$,
see also
\cite[Proposition~3.11/3.12]{Schrempf2017c9}% IEJA001 1306-6048
\ (or \cite[Proposition~2.6/2.7]{Schrempf2018a2}% X180310627 arxiv
).
This yields a natural correspondence between
factorizations and upper right zero block structure in the
system matrix (assuming zero entries in the corresponding
components of the right hand side).

\medskip
In other words: One can find (non-trivial) factors of a polynomial
by looking for ``appropriate'' transformations (of a \emph{minimal} ALS).
This is the main topic in Section~\ref{sec:ff.polyfact}
respectively 
\cite[Section~2]{Schrempf2017b9}% JSC2018 0747-7171
.
If one factorizes a polynomial in two (not necessarily irreducible)
factors, ``their'' admissible linear systems are \emph{minimal}.
The converse ---and that is the core of Section~\ref{sec:ff.polymul}---
is also true. Although the minimal polynomial multiplication
(Proposition~\ref{pro:ff.minmul})
seems to be obvious, the proof is \emph{highly} non-trivial.
(A possible reason is that only minimality is assumed and not,
for example, invertibility of the system matrix over the formal
power series.)

\cite[Example~50]{Schrempf2017b9}% JSC2018 0747-7171
\ could serve as an appetizer. There the polynomial factorization
is used to compute the eigenvalues of a matrix via the factorization
of its characteristic polynomial.

\subsection{Preliminaries}\label{sec:ff.fprelim}

For the main definitions we refer to
\cite[Section~1, Page~5]{Schrempf2017b9}% JSC2018 0747-7171
. The free associative algebra $\freeALG{\field{K}}{X}$
is a ``similarity'' \emph{unique factorization domain} (UFD).
For example, the polynomials $p=1-xy$ and $q=1-yx$ are \emph{similar}
because there exist $\tilde{p},\tilde{q} \in \freeALG{\field{K}}{X}$
such that $p\tilde{p} = \tilde{q} q$ with $p,\tilde{q}$ \emph{left coprime}
and $\tilde{p},q$ \emph{right coprime}
\cite{Cohn1963b}% MR0155851 0002-9947
.

\subsection{Minimal Polynomial Multiplication}\label{sec:ff.polymul}

As an introduction one could take
the multiplication of $x$ and $1-yx$
using Proposition~\ref{pro:ff.mul1},
see also 
\cite[Example~30]{Schrempf2017b9}% JSC2018 0747-7171
. The following lemma is needed in Section~\ref{sec:ff.minpoly}
and (the proof of) the following proposition.

\begin{lemma}[\protect{%
\cite[Lemma~3.15]{Schrempf2017c9}% IEJA001 1306-6048
}]\label{lem:ff.min1}
Let $\als{A} = (u,A,v) = (1,A,\lambda)$ be an ALS of dimension $n\ge 2$
and $\field{K}$-linearly dependent left family $s=A\inv v$.
Let $m \in \{ 2, 3, \ldots, n \}$ be the minimal index
such that the left subfamily $\underline{s} = (A\inv v)_{i=m}^n$
is $\field{K}$-linearly independent.
Let $A = (a_{ij})$ and assume that $a_{ii}=1$ for $1 \le i \le m$
and $a_{ij}=0$ for $j < i \le m$ (upper triangular $m \times m$ block)
and $a_{ij}=0$ for $j \le m < i$ (lower left zero block of size $(n-m) \times m$).
Then there exists matrices $T,U \in \field{K}^{1 \times (n+1-m)}$
such that
\begin{displaymath}
U + (a_{m-1,j})_{j=m}^n - T(a_{ij})_{i,j=m}^n  =
\begin{bmatrix}
  0 & \ldots & 0
\end{bmatrix}
\quad\text{and}\quad
T(v_i)_{i=m}^n = 0.
\end{displaymath}
\end{lemma}

\begin{proposition}[Minimal Polynomial Multiplication
\protect{\cite[Proposition~28]{Schrempf2017b9}% JSC2018 0747-7171
}]\label{pro:ff.minmul}
Let $p,q \in \freeALG{\field{K}}{X}\setminus \field{K}$ be given by the
\emph{minimal} polynomial admissible linear systems
$A_p = (1, A_p, \lambda_p)$ of dimension $n_f$ and
$A_q = (1, A_q, \lambda_q)$ of dimension $n_g$ respectively.
Then the ALS $\als{A}$ from Proposition~\ref{pro:ff.mul1} for $pq$
is \emph{minimal} of dimension $n = n_p + n_q - 1$.
\end{proposition}

\subsection{Polynomial Factorization}\label{sec:ff.polyfact}

The polynomial factorization theory depends on \emph{minimal}
(polynomial) admissible linear systems. How to obtain such
systems directly is discussed in Section~\ref{sec:ff.minsys}.
How to construct them in general is discussed in
Section~\ref{sec:ff.minpoly}.

\begin{remark}
Notice that, although we use (general) admissible linear systems
here to represent polynomials, the factorization does \emph{not}
depend on the construction of the free field. Indeed, the
system matrix of a \emph{minimal} linear representation
of a polynomial is already invertible over the free associative
algebra.
\end{remark}

\begin{theorem}[Polynomial Factorization
\protect{\cite[Theorem~40]{Schrempf2017b9}% JSC2018 0747-7171
}]\label{thr:ff.factorization}
Let $p \in \freeALG{\field{K}}{X}$ be given by the
\emph{minimal} polynomial admissible linear system $\als{A} = (1,A,\lambda)$
of dimension $n = \rank p \ge 3$.
Then $p$ factorizes in $p = q_1 q_2$ with $\rank(q_i) = n_i \ge 2$
if and only if there exists a polynomial factorization transformation $(P,Q)$
such that $PAQ$ has an upper right block of zeros of size $(n_1 - 1) \times (n_2 - 1)$.
\end{theorem}

As we have already seen in the beginning of this section,
we have to find an admissible transformation (over the ground
field $\field{K}$) to create
upper right blocks of zeros (of appropriate size) in a
minimal polynomial ALS to detect (non-trivial) factors
of a polynomial.
If $\field{K}$ is not algebraically closed it can be difficult
to check if there is a solution. A simple case is illustrated
in \cite[Example~37]{Schrempf2017b9}% JSC2018 0747-7171
.
The practical application is by
\cite[Proposition~42]{Schrempf2017b9}% JSC2018 0747-7171
, a simple variant of 
\cite[Theorem~4.1]{Cohn1999a}% MR1723470 0218-1967
.

\subsection{Factorization Theory$^{\star\star}$}\label{sec:ff.ft}

The general factorization theory is somewhat difficult.
Although it seems to be clear from
the polynomials how it should be, the path to 
the \emph{divisibility equivalence}
(Theorem~\ref{thr:ff.ldivs}) is long and stony.
One needs a notion of left (respectively right) divisibility
on the level of \emph{minimal} admissible linear systems.
This is not straight forward (for details we refer to
\cite[Section~4]{Schrempf2017c9}% IEJA001 1306-6048
).
But in return one can ``forget'' the free associative algebra
and factorize elements directly in the free field.
And also here there are two sides of one coin, namely
the (minimal) multiplication in Section~\ref{sec:ff.factmul}
respectively 
\cite[Theorem~5.2]{Schrempf2017c9}% IEJA001 1306-6048
\ and the factorization via detecting zero blocks
in Section~\ref{sec:ff.genfact}
respectively
\cite[Theorem~5.9]{Schrempf2017c9}% IEJA001 1306-6048
.

\begin{theorem}[Divisibility Equivalence \protect{%
\cite[Theorem~4.11]{Schrempf2017c9}% IEJA001 1306-6048
}]\label{thr:ff.ldivs}
Let $p,q \in \freeALG{\field{K}}{X}$. Then $p$ left (respectively right)
divides $q$ if and only if $p$ left (respectively right)
divides $q$ in $\field{F} = \freeFLD{\field{K}}{X}$.
\end{theorem}

\subsection{Minimal Factor Multiplication$^{\star\star}$}\label{sec:ff.factmul}

Given two minimal admissible systems,
under which conditions are the multiplications
from Proposition~\ref{pro:ff.ratop} and~\ref{pro:ff.mul1}
\emph{minimal}? A special case is the minimal polynomial
multiplication (Proposition~\ref{pro:ff.minmul}).
The general answer is given in
\cite[Theorem~5.2]{Schrempf2017c9}% IEJA001 1306-6048
\ within the (framework of the) general factorization theory.

\subsection{General Factorization$^{\star\star}$}\label{sec:ff.genfact}

Like in the general (minimal) multiplication in the previous subsection
we have to distinguish several cases for the
factorization 
\cite[Theorem~5.9]{Schrempf2017c9}% IEJA001 1306-6048
. Looking for zero (lower left and upper right) blocks
(of appropriate size)
in the system matrix of a \emph{minimal} ALS
(similar to the polynomial factorization) is rather
natural when we want to ``reverse'' the multiplication.
The main difficulties however are far from obvious
and therefore one of the first steps in the general
factorization theory
\cite[Section~4]{Schrempf2017c9}% IEJA001 1306-6048
\ is to define, \emph{what} we mean by a ``factor''
(since in a field there are no non-zero non-units,
that is, \emph{each} non-zero element is invertible).

\subsection{Examples Factorization}\label{sec:ff.examples}

Polynomial factorization is illustrated in detail (step by step) in
\cite[Section~4]{Schrempf2017b9}% JSC2018 0747-7171
.
The general factorization (of a regular element)
is discussed briefly in 
\cite[Example~5.10]{Schrempf2017c9}% IEJA001 1306-6048
.

\section{Minimizing}\label{sec:ff.m}

The basic idea of the minimization (of a linear representation)
with left and right minimization steps is surprisingly simple.
If the block structure becomes coarser and a ``look'' is not
sufficient any more, row and column transformations can be
found by solving a \emph{linear} system of equations.
That is the essential content of Section~\ref{sec:ff.wp}
(word problem), the foundation stone of the whole theory.
The naive idea was to solve ``local'' word problems,
producing plenty of questions which ---among other things---
led to the factorization theory \ldots

But \emph{when}, that is, under which conditions,
is an admissible linear system (constructed out of two
\emph{minimal} ones by Proposition~\ref{pro:ff.ratop})
minimal? If there are no more left or right
``linear'' minimization steps possible?
Is it sufficient to find \emph{one} ``finest'' structure
such that the system matrix is an upper block triangular
matrix with a maximal number of (quadratic) diagonal blocks?

For polynomials (given by polynomial admissible linear systems)
this can be done by a relatively simple algorithm
which is formulated in Section~\ref{sec:ff.minpoly}.
If one knows ``all'' factorizations of a polynomial,
one also knows all ``finest'' pivot block structures of
the \emph{minimal} admissible linear systems of its inverse
and one can continue to calculate ``easily'' because
it is still rather simple to minimize.

Already in the beginning of Section~\ref{sec:ff.r}
(calculating) we have discussed assumptions on
the construction of an ALS for the inverse of an element.
In Section~\ref{sec:ff.verfeinern} 
we investigate the connection between a factorization
and the refinement of pivot blocks in the system of the
inverse a little more thoroughly and describe the
approach of the latter.
One of the central question in Section~\ref{sec:ff.minallg}
is that of a \emph{sufficient} condition for
the minimization with \emph{linear} techniques.

In fact one could develop a general
minimization algorithm using polynomial systems of equations.
However, these are usually difficult to solve.
And if we do not know anything about the existence of a
solution, we do not know anything about minimality.
Therefore non-linear techniques should be avoided
whenever this is possible by ``keeping'' a fine block structure.

\medskip
Since the main goal of this section is to ``minimize''
addition and multiplication, some thoughts from this
point of view should be summarized.
That the factorization of an element does make sense
for the multiplication is immediately clear:
In this case one can cancel factors.
This is used for example to find the left greatest
common divisor of two polynomials
\cite[Example~5.4]{Schrempf2018a2}% X180310627 arxiv
.
But it is not that trivial since an atom might
not necessarily lie ``beside'' its inverse,
for example
\begin{displaymath}
x (1-yx) \cdot x\inv = (1-xy)x\cdot x\inv = 1-xy.
\end{displaymath}
Additionally it can happen that two irreducible
elements ``fuse'' to one
\cite[Section~4]{Schrempf2017c9}% IEJA001 1306-6048
\ and therefore we need a refinement of
pivot blocks ``inside'' an atom (irreducible element).
But also from an additive point of view the
factorization plays a crucial role because one
needs ``common'' left and right factors of two summands
only ``once''.
Notice that there are also linear techniques for
refinement, for example to bring an ALS to a suitable
form for the minimal inverse (Theorem~\ref{thr:ff.mininv}).

\medskip
Recall that here we operate \emph{directly} in the
(system matrix of the) linear representation and
therefore we are \emph{independent} of its regularity
(that is, invertibility over the formal power series).
And that has its price.
The ``classical'' methods for the minimization of linear
representations for \emph{regular} elements work
mainly \emph{indirectly} by computing the left
and right families,
see for example
\cite[Section~3]{Schrempf2017a9}% IJAC2018 0218-1967
.

\subsection{Preliminaries and a Standard Form}

To be able to formulate statements ---in particular
for the minimization--- in a convenient way, we need
some notation which formalizes what we have already
used, namely to describe an ALS (and admissible
transformations) in terms of \emph{block} rows and
columns instead of (single) rows and columns.
Then it is possible to define a \emph{standard form}
which plays an important role when we want to minimize
admissible linear systems coming from addition or
multiplication (later in Section~\ref{sec:ff.minallg}).
This is the first part in 
\cite[Section~3]{Schrempf2018a2}% X180310627 arxiv
. To construct a \emph{standard admissible linear system}
out of a \emph{minimal} ALS we need to ``refine'' it.
This is the goal of Section~\ref{sec:ff.verfeinern},
the second part in
\cite[Section~3]{Schrempf2018a2}% X180310627 arxiv
.

\begin{definition}[Pivot Blocks, Pivot Block Transformation
\protect{\cite[Definition~3.1]{Schrempf2018a2}% X180310627 arxiv
}]\label{def:ff.pivot}
Let $\als{A} = (u,A,v)$ be an ALS
and denote $A = (A_{ij})_{i,j=1}^m$ the block decomposition
(with square diagonal blocks $A_{ii}$) with \emph{maximal} $m$
such that $A_{ij}=0$ for $i > j$.
The diagonal blocks $A_{ii}$ are called \emph{pivot blocks},
the number $m$ is denoted by $\pivot\als{A}$.
The \emph{dimension} (or \emph{size})
of a pivot block $A_{ii}$ for $i \in \{ 1, 2, \ldots, m \}$
is $n_i = \dim_i \als{A}$.
For a pivot block $k$ let $I_{1:k-1}$ (respectively $I_{k+1:m}$)
denote the identity matrix of size $n_1 + \ldots + n_{k-1}$
(respectively $n_{k+1} + \ldots + n_m$).
An admissible transformation $(P,Q)$ of the form
\begin{displaymath}
(P,Q)_k = 
\left(
\begin{bmatrix}
I_{1:k-1} & . & . \\
. & \bar{T} & . \\
. & . & I_{k+1:m}
\end{bmatrix},
\begin{bmatrix}
I_{1:k-1} & . & . \\
. & \bar{U} & . \\
. & . & I_{k+1:m}
\end{bmatrix}
\right)
\end{displaymath}
with $\bar{T},\bar{U} \in \field{K}^{n_k \times n_k}$ is called
(admissible) \emph{$k$-th pivot block transformation}.
\end{definition}

\begin{definition}[Refined Pivot Block and Refined ALS
\protect{\cite[Definition~3.3]{Schrempf2018a2}% X180310627 arxiv
}]\label{def:ff.redpivot}
Let $\als{A} = (u,A,v)$ be an ALS with $m = \pivot\als{A}$
pivot blocks of size $n_i = \dim_i\als{A}$.
A pivot block $A_{kk}$ (for $1 \le k \le m$)
is called \emph{refined}
if there does not exist an admissible pivot block transformation
$(P,Q)_k$ such that $(PAQ)_{kk}$ has a lower left
block of zeros of size $i \times (n_k -i)$
for an $i \in \{1, 2, \ldots, n_k-1\}$.
The admissible linear system $\als{A}$ is called \emph{refined}
if all pivot blocks are refined.
\end{definition}

\begin{definition}[Standard Admissible Linear System
\protect{\cite[Definition~3.8]{Schrempf2018a2}% X180310627 arxiv
}]\label{def:ff.stdals}
A \emph{minimal} and \emph{refined}
ALS $\als{A} = (u,A,v) = (1,A,\lambda)$,
that is, $v = [0,\ldots, 0, \lambda]$, is called
\emph{standard}.
\end{definition}

\begin{remark}
For a polynomial $p$ given by a \emph{standard} ALS
$\als{A}$ (of dimension $n \ge 2$) the minimal inverse
of $\als{A}$ (of dimension $n-1$) is refined if and
only if $\als{A}$ is obtained by the minimal polynomial
multiplication of its irreducible factors $q_i$ in
$p = q_1 q_2 \cdots q_m$.
For a detailed discussion of polynomial factorization
(in free associative algebras) we refer to
\cite{Schrempf2017b9}% JSC2018 0747-7171
.
\end{remark}

\subsection{The Word Problem$^\star$}\label{sec:ff.wp}

One of the difficulties in free fields is (that of) the \emph{word problem},
that is, to check whether two admissible linear systems represent
the \emph{same} element. A solution to the word problem is
\cite[Theorem~4.1]{Cohn1999a}% MR1723470 0218-1967
. Unfortunately it is hard to apply practically already
for systems of dimension~3.
If those systems are given by \emph{minimal} admissible linear systems
however, the word problem can be ``linearized'', that is,
it is equivalent to the solution of a \emph{linear} system of
equations. For a detailed discussion we refer to
\cite[Section~2]{Schrempf2017a9}% IJAC2018 0218-1967
.

\begin{theorem}[Linearized Word Problem
\protect{\cite[Theorem~2.4]{Schrempf2017a9}% IJAC2018 0218-1967
}]\label{thr:ff.wp}
Let $f,g \in \field{F}$ be given by the
\emph{minimal} admissible linear systems
$\als{A}_f = (u_f, A_f, v_f)$ and $\als{A}_g = (u_g, A_g, v_g)$
of dimension $n$ respectively.
Then $f = g$ if and only if there exist matrices $T,U \in \field{K}^{n\times n}$
such that $u_f U = 0$, $T A_g - A_f U= A_f u_f\trp u_g$ and $T v_g = v_f$.
\end{theorem}

The techniques used for the minimization in Section~\ref{sec:ff.minpoly}
and~\ref{sec:ff.minallg} can be interpreted as solving ``local''
word problems. The other way around one can view the word problem
as one ``big'' minimization step.

\subsection{Minimizing a Polynomial ALS}\label{sec:ff.minpoly}

To illustrate the main idea we (partially) minimize
a \emph{non-minimal} ``almost'' polynomial
ALS $\als{A} = (u,A,v)$ of dimension $n=6$
for $p = -xy + (xy + z)$
\cite[Section~2.2]{Schrempf2017b9}% JSC2018 0747-7171
. Notice that we do not need knowledge
of the left and right family at all.
Let
\begin{equation}\label{eqn:ff.min0}
\als{A} = \left(
\begin{bmatrix}
1 & . & . & . & . & . 
\end{bmatrix},
\begin{bmatrix}
1 & -x & . & -1 & . & . \\
. & 1 & -y & . & . & . \\
. & . &  1 & . & . & . \\
. & . & . & 1 & -x & -z \\
. & . & . & . & 1 & -y \\
. & . & . & . & . & 1 
\end{bmatrix},
\begin{bmatrix}
. \\ . \\ -1 \\ . \\ . \\ 1
\end{bmatrix}
\right).
\end{equation}
First we do one ``left'' minimization step,
that is, we remove (if possible)
one element of the $\field{K}$-linearly
dependent left family $s = A\inv v$
and construct a new system.
We fix a $1 \le k < n$, say $k=3$.
If we find a transformation $(P,Q)$ of the form
\begin{equation}\label{eqn:ff.ltrf}
(P,Q) = \left(
\begin{bmatrix}
I_{k-1} & . & . \\
. & 1 & T \\
. & . & I_{n-k}
\end{bmatrix},
\begin{bmatrix}
I_{k-1} & . & . \\
. & 1 & U \\
. & . & I_{n-k}
\end{bmatrix}
\right)
\end{equation}
such that row~$k$ in $PAQ$ is $[0,0,1,0,0,0]$
and $(Pv)_k = 0$, we can eliminate row~$k$ and
column~$k$ in $P \als{A} Q$ because $(Q\inv s)_k = 0$.
(This is what we have already done in Section~\ref{sec:ff.e}.)
How can we find these blocks $T,U \in \field{K}^{1 \times (n-k)}$?
We write $\als{A}$ in block form with respect to (block) row/column~$k$
and write $A_{1:,:m}$ for $A_{1:k-1,k+1:m}$, etc. (Recall that here
we have $m=n$ pivot blocks. Block indices are underlined to distinguish
them from component indices.)
\begin{equation}\label{eqn:ff.bloals}
\als{A}^{[\block{k}]} = \left(
\begin{bmatrix}
u_{\block{1:}} & . & . 
\end{bmatrix},
\begin{bmatrix}
A_{1:,1:} & A_{1:,k} & A_{1:,:m} \\
. & 1 & A_{k,:m} \\
. & . & A_{:m,:m}
\end{bmatrix},
\begin{bmatrix}
v_{\block{1:}} \\ v_{\block{k}} \\ v_{\block{:m}}
\end{bmatrix}
\right)
\end{equation}
and apply the transformation $(P,Q)$:
\begin{align*}
P A Q &= 
\begin{bmatrix}
I_{1:k-1} & . & . \\
. & 1 & T \\
. & . & I_{k+1:m}
\end{bmatrix}
\begin{bmatrix}
A_{1:,1:} & A_{1:,k} & A_{1:,:m} \\
. & 1 & A_{k,:m} \\
. & . & A_{:m,:m}
\end{bmatrix}
\begin{bmatrix}
I_{1:k-1} & . & . \\
. & 1 & U \\
. & . & I_{k+1:m}
\end{bmatrix} \\
&=
\begin{bmatrix}
A_{1:,1:} & A_{1:,k} & A_{1:,k} U + A_{1:,:m} \\
. & 1 & U + A_{k,:m} + T A_{:m,:m} \\
. & . & A_{:m,:m}
\end{bmatrix}, \\
Pv &= 
\begin{bmatrix}
I_{1:k-1} & . & . \\
. & 1 & T \\
. & . & I_{k+1:m}
\end{bmatrix}
\begin{bmatrix}
v_{\block{1:}} \\ v_{\block{k}} \\ v_{\block{:m}}
\end{bmatrix}
=
\begin{bmatrix}
v_{\block{1:}} \\ v_{\block{k}} + T v_{\block{:m}} \\ v_{\block{:m}}
\end{bmatrix}.
\end{align*}
Now we can read of a \emph{sufficient} condition
for $(Q\inv s)_k = 0$, namely
the \emph{existence} of $T,U \in \field{K}^{1 \times (n-k)}$ such that
\begin{equation}\label{eqn:ff.lmsys}
U + A_{k,:m} + T A_{:m,:m} = 0
\quad\text{and}\quad
v_{\block{k}} + T v_{\block{:m}} = 0.
\end{equation}
(Compare with the word problem, Theorem~\ref{thr:ff.wp}.)
Let $d$ be the number of letters in our alphabet $X$.
The blocks $T= [\alpha_{k+1}, \alpha_{k+2}, \ldots, \alpha_n]$ and
$U = [\beta_{k+1}, \beta_{k+2}, \ldots, \beta_n]$ in the transformation $(P,Q)$
are of size $1 \times (n-k)$, thus we have
a \emph{linear} system of equations (over $\field{K}$) with
$2(n-k)$ unknowns (for $k>1$) and $(d+1)(n-k) + 1$ equations:
\begin{align*}
\begin{bmatrix}
\beta_{k+1} & \beta_{k+2} & \beta_{k+3}
\end{bmatrix}
+
\begin{bmatrix}
0 & 0 & 0 
\end{bmatrix}
+ \qquad\qquad\qquad\qquad & \\
\begin{bmatrix}
\alpha_{k+1} & \alpha_{k+2} & \alpha_{k+3}
\end{bmatrix}
\begin{bmatrix}
1 & -x & -z \\
. & 1 & -y \\
. & . & 1
\end{bmatrix}
&=
\begin{bmatrix}
0 & 0 & 0
\end{bmatrix}, \\
\begin{bmatrix}
-1
\end{bmatrix}
+
\begin{bmatrix}
\alpha_{k+1} & \alpha_{k+2} & \alpha_{k+3}
\end{bmatrix}
\begin{bmatrix}
. \\ . \\ 1
\end{bmatrix}
&=
\begin{bmatrix}
0
\end{bmatrix}.
\end{align*}
One solution is $T=[0,0,1]$ and $U=[0,0,-1]$.
We compute $\tilde{\als{A}}_1 = P\als{A}Q$ 
and remove block row~$\block{k}$
and column~$\block{k}$ to get the new ALS
\begin{displaymath}
\als{A}_1 = (u,A,v) = \left(
\begin{bmatrix}
1 & . & . & . & .
\end{bmatrix},
\begin{bmatrix}
1 & -x & -1 & . & . \\
. & 1 & . & . & y \\
. & . & 1 & -x & -z \\
. & . & . & 1 & -y \\
. & . & . & . & 1
\end{bmatrix},
\begin{bmatrix}
. \\ . \\ . \\ . \\ 1
\end{bmatrix}
\right).
\end{displaymath}
For a ``right'' minimization step, that is,
removing (if possible) one element of the $\field{K}$-linearly dependent
right family $t = u A\inv$ we are looking for a transformation $(P,Q)$ of the
form
\begin{equation}\label{eqn:ff.rtrf}
(P,Q) = \left(
\begin{bmatrix}
I_{k-1} & T & . \\
. & 1 & . \\
. & . & I_{n-k}
\end{bmatrix},
\begin{bmatrix}
I_{k-1} & U & . \\
. & 1 & . \\
. & . & I_{n-k}
\end{bmatrix}
\right)
\end{equation}
such that column~$k$ in $PAQ$ is $[0,\ldots,0,1,0,\ldots,0]\trp$.
A sufficient condition for $(t P\inv)_k=0$ is the existence of
$T,U \in \field{K}^{(k-1) \times 1}$ such that
\begin{equation}\label{eqn:ff.rmsys}
A_{1:,1:} U + A_{1:,k} + T = 0.
\end{equation}
For the illustration we refer to 
\cite[Section~2.2]{Schrempf2017b9}% JSC2018 0747-7171
.
If a left (respectively right) minimization step with $k=1$
(respectively $k=n$ and $v = [0,\ldots,0,\lambda]\trp$)
can be done, then the ALS represents zero and we can stop immediately.

The following is the only non-trivial observation:
Recall that, if there exist row (respectively column) blocks $T,U$ such that
\eqref{eqn:ff.lmsys} (respectively \eqref{eqn:ff.rmsys})
has a solution then the left (respectively right) family
is $\field{K}$-linearly dependent.
To guarantee \emph{minimality} 
we need the other implication, that is,
the existence of appropriate row \emph{or} column blocks
for non-minimal polynomial admissible linear systems.

The following arguments can be found in the proof of
\cite[Proposition~28]{Schrempf2017b9}% JSC2018 0747-7171
:
Let $\als{A} = (u,A,v)$ be a polynomial ALS of
dimension $n \ge 2$ with left family $s = (s_1, s_2, \ldots, s_n)$
and assume that there exists a $1 \le k < n$ such that
the subfamily $(s_{k+1}, s_{k+2}, \ldots, s_n)$ is $\field{K}$-linearly
independent while $(s_k, s_{k+1}, \ldots, s_n)$ is
$\field{K}$-linearly dependent.
Then, by Lemma~\ref{lem:ff.min1}, there exist
matrices $T,U \in \field{K}^{1 \times (n-k)}$
such that \eqref{eqn:ff.lmsys} holds.
In other words: We have to start with $k_s=n-1$ for a left
and $k_t=2$ for a right minimization step.

If we apply one minimization step, we must check
the other family ``again'', illustrated in the following
example:
\begin{displaymath}
\als{A} = (u,A,v) = \left(
\begin{bmatrix}
1 & . & . & . & . 
\end{bmatrix},
\begin{bmatrix}
1 & -x & -y & x+y & . \\
. & 1 & . & . & -z \\
. & . & 1 & . & -z \\
. & . & . & 1 & -y \\
. & . & . & . & 1 \\
\end{bmatrix},
\begin{bmatrix}
. \\ . \\ . \\ . \\ 1
\end{bmatrix}
\right).
\end{displaymath}
Clearly, the left subfamily
$(s_3,s_4,s_5)$ 
and the right subfamily $(t_1, t_2, t_3)$
of $\als{A}$ 
are $\field{K}$-linearly independent respectively.
If we subtract row~3 from row~2 and add column~2 to column~3,
we get the ALS
\begin{displaymath}
\als{A}' = (u',A',v') = \left(
\begin{bmatrix}
1 & . & . & . & . 
\end{bmatrix},
\begin{bmatrix}
1 & -x & -x-y & x+y & . \\
. & 1 & 0 & . & 0 \\
. & . & 1 & . & -z \\
. & . & . & 1 & -y \\
. & . & . & . & 1 \\
\end{bmatrix},
\begin{bmatrix}
. \\ . \\ . \\ . \\ 1
\end{bmatrix}
\right).
\end{displaymath}
The right subfamily $(t_1'',t_2'',t_3'')$
of $\als{A}'' = \als{A}'\mthstrut^{[-2]}$
is (here) \emph{not} $\field{K}$-linearly independent anymore,
therefore we must check for a right minimization step
for $k=3$ again.

\medskip
\begin{definition}[Minimization Equations
\protect{\cite[Definition~31]{Schrempf2017b9}% JSC2018 0747-7171
}]
Let $\als{A} = (u,A,v)$ be a polynomial ALS of dimension $n \ge 2$.
Recalling the block decomposition \eqref{eqn:ff.bloals},
we denoty by $\als{A}^{[-\block{k}]}$ the ALS $\als{A}^{[\block{k}]}$
without (block) row/column~$k$ (of dimension $n-1$):
\begin{displaymath}
\als{A}^{[-\block{k}]} = \left(
\begin{bmatrix}
u_{\block{1:}} & . 
\end{bmatrix},
\begin{bmatrix}
A_{1:,1:} & A_{1:,:n} \\
. & A_{:n,:n}
\end{bmatrix},
\begin{bmatrix}
v_{\block{1:}} \\ v_{\block{:n}}
\end{bmatrix}
\right).
\end{displaymath}
For $k \in \{ 1,2,\ldots, n-1 \}$ the equations 
$U + A_{k,:n} + T A_{:n,:n} = 0$ and $v_{\block{k}} + T v_{\block{:n}} = 0$,
see \eqref{eqn:ff.lmsys},
with respect to the block decomposition $\als{A}^{[\block{k}]}$
are called \emph{left minimization equations},
denoted by $\mathcal{L}_k = \mathcal{L}_k(\als{A})$.
A solution by the row block pair $(T,U)$ is denoted by
$\mathcal{L}_k(T,U) = 0$, 
the corresponding transformation by $\bigl(P(T), Q(U) \bigr)$.
For $k \in \{ 2,3,\ldots, n \}$ the equations
$A_{1:,1:} U + A_{1:,k} + T = 0$, see \eqref{eqn:ff.rmsys},
with respect to the block decomposition $\als{A}^{[\block{k}]}$
are called \emph{right minimization equations},
denoted by $\mathcal{R}_k = \mathcal{R}_k(\als{A})$.
A solution by the column block pair $(T,U)$ is denoted by
$\mathcal{R}_k(T,U) = 0$,
the corresponding transformation by $\bigl(P(T), Q(U) \bigr)$.
\end{definition}

\begin{algorithm}[Minimizing a polynomial ALS
\protect{\cite[Algorithm~32]{Schrempf2017b9}% JSC2018 0747-7171
}]\label{alg:ff.minals}
\ \\
Input: $\als{A} = (u,A,v)$ polynomial ALS
  of dimension $n \ge 2$ (for some polynomial $p$).\\
Output: $\als{A}' = (,,)$ if $p=0$ or
        a minimal polynomial ALS $\als{A}' = (u',A',v')$ if $p \neq 0$.
        
\begin{algtest}
\hbox{}\\[-3ex]
\lnum{1:}\>$k := 2$ \\
\lnum{2:}\>while $k \le \dim \als{A}$ do \\
\lnum{3:}\>\>$n := \dim(\als{A})$ \\
\lnum{4:}\>\>$k' := n  +1 - k$ \\
\lnum{  }\>\>\textnormal{Is the left subfamily
  \raisebox{0pt}[0pt][0pt]{%
    $(s_{k'},\overbrace{\mthstrut s_{k'+1},\ldots, s_{n}}^{\text{lin.~indep.}})$}
    $\field{K}$-linearly dependent?} \\
\lnum{5:}\>\>if $\exists\, T,U \in \field{K}^{1 \times (k-1)}
  \textnormal{ admissible}: \mathcal{L}_{k'}(\als{A})=\mathcal{L}_{k'}(T,U)=0$ then \\
\lnum{6:}\>\>\>if $k' = 1$ then \\
\lnum{7:}\>\>\>\>return $(,,)$ \\
\lnum{  }\>\>\>endif \\
\lnum{8:}\>\>\>\raisebox{0pt}[0pt][0pt]{%
      $\als{A} := \bigl(P(T) \als{A} Q(U)\bigr)\mthstrut^{[-k']}$} \\
\lnum{9:}\>\>\>if $k > \max \bigl\{ 2, \frac{n+1}{2} \bigr\}$ then \\
\lnum{10:}\>\>\>\>$k := k-1$ \\
\lnum{   }\>\>\>endif \\
\lnum{11:}\>\>\>continue \\
\lnum{   }\>\>endif \\
\lnum{12--15:}\>\>\textnormal{(not necessary here)}\\
\lnum{   }\>\>\textnormal{Is the right subfamily
      \raisebox{0pt}[0pt][0pt]{%
            $(\overbrace{\mthstrut t_1, \ldots,t_{k-1}}^{\text{lin.~indep.}}, t_k)$}
          $\field{K}$-linearly dependent?} \\
\lnum{16:}\>\>if $\exists\, T,U \in \field{K}^{(k-1) \times 1}
    \textnormal{ admissible} : 
     \mathcal{R}_k(\als{A}) = \mathcal{R}_k(T,U)=0$ then \\
\lnum{17:}\>\>\>$\als{A} := \bigl(P(T) \als{A} Q(U) \bigr)\mthstrut^{[-k]}$ \\
\lnum{18:}\>\>\>if $k > \max \bigl\{ 2, \frac{n+1}{2} \bigr\}$ then \\
\lnum{19:}\>\>\>\>$k := k-1$ \\
\lnum{   }\>\>\>endif \\
\lnum{20:}\>\>\>continue \\
\lnum{   }\>\>endif \\
\lnum{21:}\>\>$k := k+1$ \\
\lnum{   }\>done \\
\lnum{22:}\>return $P\als{A},$ 
        \textnormal{with $P$, such that $Pv = [0,\ldots,0,\lambda]\trp$}
\end{algtest}
\end{algorithm}

\begin{remark}
Notice that, compared to 
\cite[Algorithm~4.14]{Schrempf2018a2}% X180310627 arxiv
, the first row does not have to be treated separately
(using an \emph{extended} ALS), because for $\dim \als{A} =2$
$\field{K}$-linear independence of the left family is
equivalent to $\field{K}$-linear independence of the right
family. Hence the former is indirectly checked by the
latter in line~16 and therefore the lines 12--15
(in the general algorithm) do not have a correspondence here.
\end{remark}

\subsection{Pivot Block Refinement}\label{sec:ff.verfeinern}

To be able to minimize an ALS using \emph{linear} techniques
only the pivot blocks have to be \emph{refined},
that is, none can be (admissibly) transformed such that it
splits in two (smaller) pivot blocks.
For an illustration we consider the ALS
\begin{displaymath}
\als{A} = \left(
\begin{bmatrix}
1 & . & . & . 
\end{bmatrix},
\begin{bmatrix}
1 & -z & . & . \\
. & 2+x & . & 1 \\
. & 2y & -3 & y \\
. & x & 3x & .
\end{bmatrix},
\begin{bmatrix}
. \\ . \\ . \\ 1
\end{bmatrix}
\right)
\end{displaymath}
with a $3 \times 3$ pivot block.
Using the admissible block transformation
\begin{displaymath}
(P,Q) = \left(
\begin{bmatrix}
1 & . & . & . \\
. & \alpha_{2,2} & \alpha_{2,3} & 0 \\
. & \alpha_{3,2} & \alpha_{3,3} & 0 \\
. & \alpha_{4,2} & \alpha_{4,3} & 1
\end{bmatrix},
\begin{bmatrix}
1 & . & . & . \\
. & \beta_{2,2} & \beta_{2,3} & \beta_{2,4} \\
. & \beta_{3,2} & \beta_{3,3} & \beta_{3,4} \\
. & \beta_{4,2} & \beta_{4,3} & \beta_{4,4} 
\end{bmatrix}
\right)
\end{displaymath}
we need to check if it is possible to create
a lower left block of zeros of size $1 \times 2$
or $2 \times 1$ in the second pivot block of $P \als{A} Q$.
First we need to ensure invertibility of $P$ and $Q$ by
the conditions
\begin{align*}
0 &\neq \det(P) = \alpha_{2,2}\alpha_{3,3} - \alpha_{2,3}\alpha_{3,2}
  \quad\text{and}\\
0 &\neq \det(Q) =
          (\beta_{2,2}\beta_{3,3} - \beta_{2,3}\beta_{3,2})\beta_{4,4} \\
        & \qquad\qquad\qquad\qquad
        + (\beta_{2,4}\beta_{3,2} - \beta_{2,2}\beta_{3,4})\beta_{4,3}
        + (\beta_{2,3}\beta_{3,4} - \beta_{2,4}\beta_{3,3})\beta_{4,2}.
\end{align*}
To (possibly) split the second pivot block into a $1\times 1$ and $2\times 2$
block we need to solve the equations obtained by applying
the block transformation matrices to the corresponding
coefficient matrices for $1$, $x$ and $y$
(notice that there is no contribution with respect to $z$;
irrelevant equations are marked with ``*'' on the right hand side)
\begin{align*}
\begin{bmatrix}
\alpha_{2,2} & \alpha_{2,3} & 0 \\
\alpha_{3,2} & \alpha_{3,3} & 0 \\
\alpha_{4,2} & \alpha_{4,3} & 1
\end{bmatrix}
\begin{bmatrix}
2 & . & 1 \\
. & -3 & . \\
. & . & . 
\end{bmatrix}
\begin{bmatrix}
\beta_{2,2} & \beta_{2,3} & \beta_{2,4} \\
\beta_{3,2} & \beta_{3,3} & \beta_{3,4} \\
\beta_{4,2} & \beta_{4,3} & \beta_{4,4} 
\end{bmatrix}
&= \begin{bmatrix}
* & * & * \\
0 & * & * \\
0 & * & *
\end{bmatrix}
\quad\text{for $1$, and}\\
\begin{bmatrix}
\alpha_{2,2} & \alpha_{2,3} & 0 \\
\alpha_{3,2} & \alpha_{3,3} & 0 \\
\alpha_{4,2} & \alpha_{4,3} & 1
\end{bmatrix}
\begin{bmatrix}
1 & . & . \\
. & . & . \\
1 & 3 & . 
\end{bmatrix}
\begin{bmatrix}
\beta_{2,2} & \beta_{2,3} & \beta_{2,4} \\
\beta_{3,2} & \beta_{3,3} & \beta_{3,4} \\
\beta_{4,2} & \beta_{4,3} & \beta_{4,4} 
\end{bmatrix}
&= \begin{bmatrix}
* & * & * \\
0 & * & * \\
0 & * & *
\end{bmatrix}
\quad\text{for $x$, and}\\
\begin{bmatrix}
\alpha_{2,2} & \alpha_{2,3} & 0 \\
\alpha_{3,2} & \alpha_{3,3} & 0 \\
\alpha_{4,2} & \alpha_{4,3} & 1
\end{bmatrix}
\begin{bmatrix}
. & . & . \\
2 & . & 1 \\
. & . & . 
\end{bmatrix}
\begin{bmatrix}
\beta_{2,2} & \beta_{2,3} & \beta_{2,4} \\
\beta_{3,2} & \beta_{3,3} & \beta_{3,4} \\
\beta_{4,2} & \beta_{4,3} & \beta_{4,4} 
\end{bmatrix}
&= \begin{bmatrix}
* & * & * \\
0 & * & * \\
0 & * & *
\end{bmatrix}
\quad\text{for $y$.}
\end{align*}
Thus, additionally to $\det(P)=1$ and $\det(Q)=1$, we get the equations
\begin{align*}
\alpha_{3,2}\beta_{4,2} - 3 \alpha_{3,3} \beta_{3,2} + 2\alpha_{3,2}\beta_{2,2} &= 0,\\
\alpha_{4,2}\beta_{4,2} - 3 \alpha_{4,3} \beta_{3,2} + 2\alpha_{4,2}\beta_{2,2} &= 0,\\
\alpha_{3,2}\beta_{2,2} &= 0, \\
3\beta_{3,2} + (\alpha_{4,2} + 1)\beta_{2,2} &= 0, \\
\alpha_{3,3}\beta_{4,2} + 2\alpha_{3,3}\beta_{2,2} &= 0 \quad\text{and} \\
\alpha_{4,3}\beta_{4,2} + 2\alpha_{4,3}\beta_{2,2} &= 0
\end{align*}
with (at least one) solution
\begin{displaymath}
(P,Q) = \left(
\begin{bmatrix}
1 & . & . & . \\
. & 1 & 0 & 0 \\
. & 0 & 1 & 0 \\
. & -1 & 0 & 1 
\end{bmatrix},
\begin{bmatrix}
1 & . & . & . \\
. & 1 & 0 & 0 \\
. & 0 & 0 & \frac{1}{3} \\
. & -2 & 1 & 0
\end{bmatrix}
\right)
\end{displaymath}
yielding the (refined) admissible linear system
\begin{displaymath}
P \als{A} Q = \left(
\begin{bmatrix}
1 & . & . & . 
\end{bmatrix},
\begin{bmatrix}
1 & -z & . & . \\
. & x & 1 & . \\
. & . & y & -1 \\
. & . & -1 & x
\end{bmatrix},
\begin{bmatrix}
. \\ . \\ . \\ 1
\end{bmatrix}
\right)
\end{displaymath}
representing $z \bigl(x\inv (1-xy)\inv \bigr)$.
Notice that here it would also be possible to create
a lower left $1\times 2$ zero block in the second
pivot block of $\als{A}$. This would correspond to
the factorization $z \bigl((1-yx)\inv x\inv \bigr)$,
while the original ALS could be interpreted as
$z (x-xyx)\inv$.

\medskip
Solving such polynomial systems of equations in general
is very difficult, especially if the ground field $\field{K}$
is not algebraically closed, that is, $\field{K} \subsetneq \aclo{\field{K}}$.
For further information we refer to
\cite[Section~4]{Schrempf2017b9}% JSC2018 0747-7171
\ and/or 
\cite[Example~5.10]{Schrempf2017c9}% IEJA001 1306-6048
.

\begin{Remark}
To ensure invertibility of the transformation matrices
$P$ and $Q$ one can use additional (commuting) variables
$P'=(\gamma_{ij})$, $Q'=(\delta_{ij})$ and equations
$PP'=I$, $QQ'=I$ instead of $\det(P)=1$, $\det(Q)=1$.
To say anything about the difference with respect to the
computation of Groebner bases, detailled investigations
would be necessary. An introduction to the necessary
concepts is \cite{Cox2015a}% MR3330490 book
.
\end{Remark}

\begin{Remark}
Since the computation of appropriate transformation matrices
for the refinement of (unrefined) pivot blocks in general is
difficult, one should try simpler techniques (before) to split
pivot blocks. If the permutation of rows and/or columns is
not successful, \emph{linear} techniques could be used by
avoiding ``overlapping'' of row and column transformations:
As an example we take the following ALS for $(x-xyx)\inv$,
\begin{displaymath}
\begin{bmatrix}
x & 1 & . \\
1 & y-1 & -1 \\
. & x & x
\end{bmatrix}
s =
\begin{bmatrix}
. \\ . \\ 1
\end{bmatrix},
\end{displaymath}
and assume that we want to create a lower left block of zeros
of size $2 \times 1$ in the system matrix. Then the ansatz
\begin{displaymath}
(P,Q) = \left(
\begin{bmatrix}
1 & . & . \\
\alpha_{2,1} & 1 & . \\
\alpha_{3,1} & . & 1
\end{bmatrix},
\begin{bmatrix}
1 & . & . \\
0 & 1 & . \\
\beta_{3,1} & . & 1
\end{bmatrix}
\right)
\end{displaymath}
yields a \emph{linear} system of equations with a solution
$\alpha_{2,1}=0$, $\alpha_{3,1}=-1$ and $\beta_{3,1}=1$.
This approach is also recommended for the factorization of
polynomials (to create upper right blocks of zeros).
\end{Remark}

\subsection{Minimizing a Refined ALS$^\star$}\label{sec:ff.minallg}

The core of the minimization is to establish the
equivalence of minimality and the non-existence of
solutions of certain \emph{linear} systems of equations.
Firstly we need to formalize what we have already done,
namely to apply (left and right) minimization steps
(as ``solutions'' to \emph{linear} systems of equations).
This is somewhat technical (to implement) but rather simple.
The other direction is difficult, namely to show that
there is always a ``linear'' minimization step as long as
the \emph{refined} admissible linear system is not
minimal. For the theoretical details we refer to
\cite[Section~4]{Schrempf2018a2}% X180310627 arxiv
.

\medskip
The basic procedure for the minimization is similar to that in
Algorithm~\ref{alg:ff.minals}. Instead of $n$ pivot blocks of
size $1 \times 1$ (with entry~$1$) we operate with respect to
$m \le n$ (general) pivot blocks of size $n_i\times n_i$ with
$n_1 + n_2 + \ldots + n_m = n$. To illustrate the setup of a
linear system of equations to (possibly) eliminate a block
we take (again) the ALS $\als{A}=(u,A,v)$ from Example~\ref{ex:ff.blockmin}
for $ff\inv=1$ with $f=xy-z$, namely
\begin{displaymath}
\begin{bmatrix}
1 & -x & z & . & . \\
. & 1 & -y & . & . \\
. & . & 1 & -1 & . \\
. & . & . & y & -1 \\
. & . & . & -z & x
\end{bmatrix}
s =
\begin{bmatrix}
. \\ . \\ . \\ . \\ 1
\end{bmatrix}.
\end{displaymath}
However here we are minimizing in a complete systematic way.
The system matrix has $m=4$ pivot blocks of size $n_1=n_2=n_3=1$
and $n_4=2$ respectively. We start with block $k_s=m-1=3$ for
a left minimization step. Notice that the left subfamily $(s_4, s_5)$
is $\field{K}$-linearly independent because we obtained this
(sub-)system by applying the minimal inverse (Theorem~\ref{thr:ff.mininv})
on the \emph{minimal} ALS
\begin{displaymath}
\begin{bmatrix}
1 & -x & z \\
. & 1 & -y \\
. & . & 1  \\
\end{bmatrix}
s =
\begin{bmatrix}
. \\ . \\ 1
\end{bmatrix}
\end{displaymath}
for $f=xy-z$. To check if the left subfamily $(s_3, s_4, s_5)$ is
$\field{K}$-linearly independent, we look for an admissible transformation
\begin{displaymath}
(P,Q) = 
\left(
\begin{bmatrix}
1 & . & . & . & . \\
. & 1 & . & . & . \\
. & . & 1 & \alpha_{3,4} & \alpha_{3,5} \\
. & . & . & 1 & . \\
. & . & . & . & 1
\end{bmatrix},
\begin{bmatrix}
1 & . & . & . & . \\
. & 1 & . & . & . \\
. & . & 1 & \beta_{3,4} & \beta_{3,5} \\
. & . & . & 1 & . \\
. & . & . & . & 1
\end{bmatrix}
\right)
\end{displaymath}
such that $P\als{A}Q$ has the form (``$*$'' denotes an arbitrary entry)
\begin{displaymath}
\begin{bmatrix}
1 & -x & z & * & * \\
. & 1 & -y & * & * \\
. & . & 1 & 0 & 0 \\
. & . & . & y & -1 \\
. & . & . & -z & x
\end{bmatrix}
s =
\begin{bmatrix}
. \\ . \\ 0 \\ . \\ 1
\end{bmatrix}
\end{displaymath}
by solving the linear system of equations
\begin{align*}
-1 + \beta_{3,4} + \alpha_{3,4}y - \alpha_{3,5} z &=0, \\
\beta_{3,5} - \alpha_{3,4} + \alpha_{3,5} x &=0 \quad\text{and} \\
\alpha_{3,5} &= 0.
\end{align*}
Notice that these are indeed $5$~equations, namely
$\beta_{3,4}-1=0$, $\beta_{3,5}-\alpha_{3,4}=0$ and $\alpha_{3,5}=0$
(for~$1$), $\alpha_{3,5} = 0$ (for~$x$), $\alpha_{3,4}=0$ (for~$y$)
and $-\alpha_{3,5}=0$ (for $z$). Since there is a solution
($\alpha_{3,4}=\alpha_{3,5}=\beta_{3,5}=0$ and $\beta_{3,4}=1$)
the left subfamily $(s_3, s_4, s_5)$ is $\field{K}$-linearly
dependent. Applying the transformation $(P,Q)$ with the appropriate
entries on $\als{A}$ and removing row~3 and column~3 yields the
ALS $\als{A}'$,
\begin{displaymath}
\begin{bmatrix}
1 & -x & z & . \\
. & 1 & -y & . \\
. & . & y & -1 \\
. & . & -z & x
\end{bmatrix}
s' =
\begin{bmatrix}
. \\ . \\ . \\ 1
\end{bmatrix}
\end{displaymath}
with $m'=3$ pivot blocks. Now $k'_s=m'-1=2$ and we check if
the (new) left subfamily $(s'_2, s'_3, s'_4)$ is $\field{K}$-linearly
independent by looking for a transformation
\begin{displaymath}
(P',Q') = 
\left(
\begin{bmatrix}
1 & . & . & . \\
. & 1 & \alpha_{2,3} & \alpha_{2,4} \\
. & . & 1 & . \\
. & . & . & 1
\end{bmatrix},
\begin{bmatrix}
1 & . & . & . \\
. & 1 & \beta_{2,3} & \beta_{2,4} \\
. & . & 1 & . \\
. & . & . & 1
\end{bmatrix}
\right)
\end{displaymath}
such that $P'\als{A}'Q'$ has the form
\begin{displaymath}
\begin{bmatrix}
1 & -x & * & * \\
. & 1 & 0 & 0 \\
. & . & y & -1 \\
. & . & -z & x
\end{bmatrix}
s' =
\begin{bmatrix}
. \\ 0 \\ . \\ 1
\end{bmatrix}.
\end{displaymath}
Since such a transformation exists ($\alpha_{2,4}=\beta_{2,3}=0$,
$\alpha_{2,3}=\beta_{2,4}=1$), the left subfamily $(s'_2, s'_3, s'_4)$
is $\field{K}$-linearly dependent. By removing row~2 and column~2 from
$P'\als{A}'Q'$ we obtain the ALS $\als{A}''= (P'\als{A}'Q')^{[-\block{k}]}$,
\begin{displaymath}
\begin{bmatrix}
1 & z & -x \\
. & y & -1 \\
. & -z & x
\end{bmatrix}
s'' =
\begin{bmatrix}
. \\ . \\ 1
\end{bmatrix}
\end{displaymath}
with $m''=2$ pivot blocks.
Now $k''_s=m''-1$ and it turns out that the (new) left
subfamily $(s''_1, s''_2, s''_3)$ is $\field{K}$-linearly
independent because there is no solution to the corresponding
linear system of equations. Therefore we switch to the
right family and try a right minimization step for $k''_t=2$,
that is, looking for a transformation
\begin{displaymath}
(P'',Q'') = \left(
\begin{bmatrix}
1 & \alpha_{1,2} & \alpha_{1,3} \\
. & 1 & . \\
. & . & 1
\end{bmatrix},
\begin{bmatrix}
1 & \beta_{1,2} & \beta_{1,3} \\
. & 1 & . \\
. & . & 1
\end{bmatrix}
\right)
\end{displaymath}
such that $P''\als{A}''Q''$ has the form
\begin{displaymath}
\begin{bmatrix}
1 & 0 & 0 \\
. & y & -1 \\
. & -z & x
\end{bmatrix}
s'' =
\begin{bmatrix}
* \\ . \\ 1
\end{bmatrix}.
\end{displaymath}
Notice that the entries $\beta_{1,j}$ in the first row of $Q''$ have to be
zero for $(P'',Q'')$ to be admissible and the corresponding entries in the left hand
side of
\begin{displaymath}
\begin{bmatrix}
1 & 0 & 0
\end{bmatrix}
= t''
\begin{bmatrix}
1 & 0 & 0 \\
. & y & -1 \\
. & -z & x
\end{bmatrix}
\end{displaymath}
are always zero. Thus the linear system of equations for checking
$\field{K}$-linear independence of the right (sub-)family
$t''=(t''_1, t''_2, t''_3)$ is
\begin{align*}
z+\alpha_{1,2}y-\alpha_{1,3}z &= 0 \quad\text{and} \\
-x-\alpha_{1,2}+\alpha_{1,3}x &= 0.
\end{align*}
Since it has a solution for $\alpha_{1,2}=0$ and $\alpha_{1,3}=1$,
the right family $t''$ is $\field{K}$-linearly dependent.
Removing block row~2 and block column~2 from $P''\als{A}''Q''$
yields the \emph{minimal} ALS $\als{A}'''=(1,[1],1)$ for 
$f f\inv=1$. Although minimality is obvious here, it is the main
result of the (general) minimization algorithm, given a
\emph{refined} admissible linear system.
For the case $f\inv f=1$ one has to treat the first block row
seperately by using an \emph{extended} ALS. This is illustrated
in \cite[Section~4]{Schrempf2018a2}% X180310627 arxiv
.

\begin{Remark}[Correction of \protect{%
\cite[Algorithm~4.5.15]{Schrempf2018b}% TH105 thesis
}]\label{rem:ff.corralg}
As the following example shows, it is neccessary (in the
general case) to decrement the counter $k$ if a left
block minimization step was successful for $k=m$ and $k>2$,
that is, to insert lines 9--10 after line~14.
Since the left family $s$ of
\begin{displaymath}
\begin{bmatrix}
1 & -1 & . & . \\
. & 1 & x & . \\
. & -y & 1 & -1 \\
. & . & . & 1
\end{bmatrix}
s =
\begin{bmatrix}
. \\ . \\ . \\ 1
\end{bmatrix}
\end{displaymath}
is $\field{K}$-linearly dependent while the
subfamily $(s_{\block{2}},s_{\block{3}}) = (s_2,s_3,s_4)$
is $\field{K}$-linearly independent, the first (block) row can
be eliminated using an extended ALS. If $k$ is not decremented
(to $k=2$) in this case, the resulting ALS would be non-minimal,
since it has only $m'=2$ pivot blocks left:
\begin{displaymath}
\begin{bmatrix}
1 & x & . \\
-y & 1 & -1 \\
. & . & 1
\end{bmatrix}
s' =
\begin{bmatrix}
. \\ . \\ 1
\end{bmatrix}.
\end{displaymath}
This corresponds to lines 15--16
in (the corrected) \cite[Algorithm~4.14]{Schrempf2018a2}% X180310627 arxiv
.
\end{Remark}

\section*{Epilogue}

Learning to compute with fractions at school takes some time
and needs hard work by \emph{hand}. This will not be different
for \emph{free fractions} (but in general much more laborious).
For those who want to experiment in computer algebra systems:
The experimental implementation FDALG ``Free Division Algebra''
building on LINPEN ``Linear Multivariate Matrix Pencil''
is available in 
\cite{FRICAS2019}% FRICAS19 manual
\ since Release~1.3.5.

\section*{Acknowledgement}

I thank Karl Auinger, Soumyashant Nayak, Bill Page and Raymond Rogers
for their respective feedback about drafts of ``free fractions''
respectively implementational details. 

\bibliographystyle{alpha}
\bibliography{doku}
\addcontentsline{toc}{section}{Bibliography}

\end{document}